\documentclass[journal]{ieeetran}

\usepackage{amsmath,amssymb,amsfonts,mathtools}
\usepackage{algorithmic}
\usepackage{subcaption}
\usepackage{graphicx}
\usepackage{algorithm,algorithmic}
\usepackage[hidelinks]{hyperref}
\usepackage{arydshln}
\usepackage{bm}
\usepackage{textcomp}
\usepackage{cases}
\usepackage[inline]{enumitem}
\usepackage{tikz,pgfplots}
\usepackage{url}
\usepackage{xcolor}
\pgfplotsset{compat=newest}
\pgfplotsset{plot coordinates/math parser=false}
\pgfplotsset{every axis plot/.append style={line width=1.5pt}}
\tikzset{every picture/.style={line width=1.2pt}}
\tikzset{>=latex}
\usetikzlibrary{calc,backgrounds}

\newlength\fheight
\newlength\fwidth

\def\arraystretchval{1.2}

\newcommand{\hop}{\mathsf{H} }
\DeclareMathOperator{\rank}{rank}
\DeclareMathOperator{\im}{im}
\DeclareMathOperator{\maximize}{max}
\DeclareMathOperator{\minimize}{min}
\DeclareMathOperator{\st}{s. t.}
\DeclareMathOperator{\diag}{diag}
\DeclareMathOperator{\trace}{tr}

\newtheorem{definition}{Definition}
\newtheorem{remark}{Remark}
\newtheorem{example}{Example}

\newtheorem{assumption}{Assumption}
\newtheorem{theorem}{Theorem}
\newtheorem{lemma}{Lemma}
\newtheorem{proposition}{Proposition}

\begin{document}
\title{From~a~Frequency-Domain~Willems'~Lemma to~Data-Driven~Predictive~Control}
\author{{Tomas}~J.~{Meijer}, \IEEEmembership{Member, IEEE}, {Koen}~J.~A.~{Scheres}, \IEEEmembership{Member, IEEE}, {Sven}~A.~N.~{Nouwens},\\
{Victor}~S.~{Dolk}, and W.~P.~M.~H.~{(Maurice)}~{Heemels}, \IEEEmembership{Fellow, IEEE}
\thanks{This research received funding from the European Research Council (ERC) under the Advanced ERC grant agreement PROACTHIS, no. 101055384.}%
\thanks{The authors are with the Control Systems Technology Section, Department of Mechanical Engineering, Eindhoven University of Technology, P.O. Box 513, 5600 MB Eindhoven, The Netherlands (e-mail: \{t.j.meijer; k.j.a.scheres; s.a.n.nouwens; v.s.dolk; m.heemels\}@tue.nl).}}

\maketitle

\begin{abstract}
Willems' fundamental lemma has recently received an impressive amount of attention from the data-driven control community. In this paper, we formulate a version of this celebrated result based on frequency-domain data. In doing so, we bridge the gap between recent developments in data-driven control, and the readily-available techniques and expertise for non-parametric frequency-domain identification. We also generalize our results to combine multiple frequency-domain data sets to form a sufficiently rich data set. Building on these results, we propose a data-driven predictive control scheme based on measured frequency-domain data of the plant. This novel scheme provides a frequency-domain counterpart of the well-known data-enabled predictive control scheme DeePC based on time-domain data. Under appropriate conditions, the new frequency-domain data-driven predictive control (FreePC) scheme is equivalent to the corresponding DeePC scheme. We demonstrate the benefits of FreePC and the use of frequency-domain data in several examples and a numerical case study, including the ability to collect data in closed loop, computational benefits, and intuitive visualization of the data. 
\end{abstract}


\section{Introduction}
\label{sec:introduction} 
\IEEEPARstart{W}{illems'} (fundamental) lemma (WFL) states that, for a linear time-invariant (LTI) system, a single data sequence, which consists of $N$ time-domain input-output data points with an input that is persistently exciting (PE), can be used to characterize all possible input-output solutions of length $L<N$ of the system~\cite{Berberich2020,Verhoek2021}. An immediate consequence of WFL, which was originally introduced in~\cite{Willems2005}, is that such a persistently-exciting data sequence describes the entire behavior of the data-generating system~\cite{vanWaarde2020}. This result does not only enable the system to be identified using subspace identification methods~\cite{Verhaegen1992}, but has also caused WFL to become one of the cornerstones of modern data-driven control. The desire to extend these remarkable results to more diverse and increasingly more complicated systems has inspired a considerable body of work aimed at generalizing/adapting WFL to accommodate such systems.  Some examples include nonlinear systems~\cite{Berberich2020,Molodchyk2024}, linear parameter-varying systems~\cite{Verhoek2021}, descriptor systems~\cite{Schmitz2022,Faulwasser2023}, stochastic systems~\cite{Faulwasser2023,Pan2022}, and continuous-time systems~\cite{Rapisarda2023}, to name a few. We refer to~\cite[Table 1]{Faulwasser2023} for a more complete overview. 

One particularly impactful application of the original time-domain WFL in recent years has been data-enabled predictive control (DeePC), which uses WFL as a data-driven substitute for the prediction model in predictive control~\cite{Coulson2019}, see~\cite{Verheijen2023} for a recent survey on this lively topic. The success of DeePC has, in turn, inspired many important extensions to accommodate, e.g., descriptor systems~\cite{Schmitz2022}, linear parameter-varying systems~\cite{Verhoek2021b}, bilinear systems~\cite{Bilgic2024}, nonlinear systems~\cite{Berberich2022b,Alsati2023,Lazar2023,Azarbahram2024}, and stochastic systems~\cite{Pan2023,Breschi2023}. 

Despite the central role of WFL in data-driven control, all existing formulations rely exclusively on time-domain data. This creates a disconnect with classical control practice, where design, tuning, and identification are often performed in the frequency domain~\cite{Franklin2010,Skogestad2005}. In many industrial applications, non-parametric frequency-response function (FRF) measurements are easy to obtain, come with well-understood uncertainty quantification tools, and can be collected in closed loop to deal with unstable systems~\cite{Soderstrom1989,Pintelon2012}. However, these rich and widely available frequency-domain data sets cannot be used directly for modern control methods or WFL-based data-driven methods. Instead, they must first be converted to, respectively, a parametric (state-space) model or long time-domain trajectories, which introduces additional modelling effort, potential bias, and loss of information. In the context of experiment design and system identification, frequency-domain techniques provide some powerful benefits over time-domain techniques, see, e.g.,~\cite{Schoukens2004}, including selection of active frequency bands, identification of unstable models, and intuitive handling and visualization of uncertainty in the data. A first work aiming to obtain a frequency-domain version of WFL is~\cite{Ferizbegovic2021}. Unfortunately, this interesting result is only shown to apply to stable systems, and it does not account for conjugate symmetry in the data that arises due to the data-generating system being real. Another relevant work connecting WFL to the frequency domain is~\cite{Markovsky2024}, in which time-domain data is used to compute the frequency response of the system through WFL. This essentially complements the problem addressed in the present paper of using frequency-domain data to characterize the time-domain behavior of the system. Some noteworthy connections between predictive control and frequency domain are found in, e.g.,~\cite{Ozkan2012,Burgos2014,Shah2013}, which present frequency-domain \emph{tuning} approaches, and~\cite{Sathyanarayanan2023} which presents a scheme based on wavelets that can extract frequency content from time-domain data. Nonetheless, generally speaking, model/data-driven predictive control are time-domain control techniques, and direct use of frequency-domain data/models of the plant (without the need to turn it into a parametric (state-space) model) is not yet possible.

In this paper, we present a frequency-domain version of WFL to bridge this gap and enable the \emph{direct} use of raw frequency-domain data in data-driven analysis and control. Some key challenges in doing so are \begin{enumerate*}[label=(\alph*)] \item to devise a suitably-structured data matrix, which contains \emph{raw} frequency-domain data and appropriately accounts for conjugate symmetry of the frequency-domain data, \item to define the notion of persistence of excitation directly in terms of raw frequency-domain data, and \item to come up with a proof that, unlike~\cite{Ferizbegovic2021}, also applies to \emph{unstable systems}. \end{enumerate*} Compared to the preliminary conference version of this work~\cite{Meijer2024-nmpc}, we not only include detailed proofs of the frequency-domain WFL, but we also present an extension that allows the use of multiple data sets, which is particularly relevant when dealing with multi-input multi-output (MIMO) systems. This extension complements a similar extension to multiple data sets presented in~\cite{vanWaarde2020} for the original time-domain WFL. We also illustrate several important applications of our results, which are new compared to the preliminary version~\cite{Meijer2024-nmpc}, including frequency-domain data-driven simulation, linear quadratic regulator (LQR) design, and a method to evaluate the transfer function of the unknown system at a desired complex frequency based on a finite amount of frequency-domain data. The latter complements a similar result based on time-domain data found in~\cite{Markovsky2024}. Finally, we use our results to propose a frequency-domain data-driven predictive control (FreePC) scheme, which provides a frequency-domain counterpart of the celebrated DeePC scheme that is based on the original time-domain version of WFL. In fact, we prove that, under appropriate conditions, the new FreePC scheme is equivalent to the corresponding DeePC scheme. We also present a numerical case study to demonstrate the key benefits of FreePC and the use of frequency-domain data, including the ability to collect data in closed loop with a pre-stabilizing controller so that our results can also be applied to unstable systems, dealing with noisy data without increasing computational complexity, and intuitively visualizing the uncertainty in the data as a result of noise and transient phenomena.

The remainder of this paper is organized as follows. In Section~\ref{sec:preliminaries}, relevant preliminaries and notations are introduced. Section~\ref{sec:problem-setting} formalizes the considered problem setting. Section~\ref{sec:frequency-domain-Willems} and Section~\ref{sec:data-driven-predictive-control} present, respectively, the frequency-domain WFL and a frequency-domain data-driven predictive control scheme, which form our main contributions. In Section~\ref{sec:applications}, we investigate several interesting applications of the frequency-domain WFL. Finally, Section~\ref{sec:conclusions} provides our conclusions, and proofs of our results can be found in the Appendix.

\section{Preliminaries}
\label{sec:preliminaries}
\subsection{Basic notation}
Let $\mathbb{R}$ denote the field of reals, $\mathbb{C}$ the complex plane, and $\mathbb{Z}$ the integers. We denote $\mathbb{Z}_{[n,m]}=\{n,n+1,\hdots,m\}$ and $\mathbb{Z}_{\geqslant n}=\{n,n+1,\hdots\}$, where $n,m\in\mathbb{Z}$ and $n\geqslant m$. The imaginary unit is denoted $j$, i.e., $j^2=-1$. For a complex-valued matrix $A\in\mathbb{C}^{n\times m}$, $A^\top$, $A^\hop$, $A^*$ and $A^+$ denote its transpose, its complex-conjugate transpose, its complex conjugate, and its pseudo-inverse, respectively, whereas $\Re A$ and $\Im A$ denote its real and imaginary part. $A_{\perp}$ denotes any matrix whose columns form a basis for $\ker A$. Let $(u_1,u_2,\hdots,u_q)\coloneqq \operatorname{col}(u_1,u_2,\hdots,u_q)$ for any vectors $u_i\in\mathbb{C}^{n_i}$, $i\in\{1,2,\hdots,q\}$. We use $\otimes$ to denote the Kronecker product, and $\bm{e}_i=(0_{i-1},1,0_{n-i})\in\mathbb{R}^{n}$, $i\in\mathbb{Z}_{[1,n]}$, to denote the $i$-th $n$-dimensional elementary basis vector. 

\subsection{Persistence of excitation}
Consider the (time-domain) signal $v~\colon~\mathbb{Z}\rightarrow\mathbb{R}^{n_v}$. We denote with $v_{[r,s]}$, $r,s\in\mathbb{Z}$ with $r\leqslant s$, the vectorized restriction of $v$ to the interval $\mathbb{Z}_{[r,s]}$, i.e.,
\begin{equation*}
    v_{[r,s]}\coloneqq \begin{bmatrix}
        v_r^\top & v_{r+1}^\top & \hdots & v_{s}^\top
    \end{bmatrix}^\top.
\end{equation*}
With some abuse of notation, we also use the notation $v_{[r,s]}$ to refer to the real-valued sequence $\{v_k\}_{k\in\mathbb{Z}_{[r,s]}}$. We define the Hankel matrix of depth $L\leqslant s-r+1$ generated by $v_{[r,s]}$ as
\begin{equation*}
    H_L(v_{[r,s]})\coloneqq \begin{bmatrix}
        v_r & v_{r+1} & \hdots & v_{s-L+1}\\
        v_{r+1} & v_{r+2} & \hdots & v_{s-L+2}\\
        \vdots & \vdots & \ddots & \vdots\\
        v_{r+L-1} & v_{r+L} & \hdots & v_{s}
    \end{bmatrix}.
\end{equation*}

Next, we recall the definition of persistence of excitation.
\begin{definition}\label{def:td-PE} 
    The real-valued sequence $v_{[0,N-1]}$ is said to be persistently exciting (PE) of order $L\in\mathbb{Z}_{[1,N]}$, if the matrix $H_L(v_{[0,N-1]})$ has full row rank.
\end{definition}

\subsection{Willems' fundamental lemma}
Consider the LTI system $\Sigma$ described by
\begin{subnumcases}{\label{eq:td-system}\Sigma~\colon~}
    x_{k+1} &= $Ax_k + Bu_k,$\label{eq:td-system-state}\\
    y_k &= $Cx_k + Du_k,$\label{eq:td-system-output}
\end{subnumcases}
where $x_k\in\mathbb{R}^{n_x}$, $u_k\in\mathbb{R}^{n_u}$ and $y_k\in\mathbb{R}^{n_y}$ denote, respectively, the state, control input, and measured output (signals) of $\Sigma$ at time $k\in\mathbb{Z}$. Let $\ell_{\Sigma}\in\mathbb{Z}_{[1,n_x]}$ denote the observability index of $\Sigma$, i.e., the smallest integer for which the observability matrix 
\begin{equation}\label{eq:obs-mat}
    \mathcal{O}_{\ell_{\Sigma}} \coloneqq \begin{bmatrix}
        C^\top & (CA)^\top & \hdots & (CA^{\ell_{\Sigma}-1})^\top
    \end{bmatrix}^\top
\end{equation}
satisfies $\rank\mathcal{O}_{\ell_{\Sigma}} = \rank\mathcal{O}_{n_x}$. The system $\Sigma$ is assumed to be controllable and unknown, i.e., the matrices $(A,B,C,D)$ are unknown.
\begin{assumption}\label{asm:ctrb}
    The pair $(A,B)$ is controllable.
\end{assumption}

As mentioned before, WFL can be used to characterize the entire behavior of $\Sigma$ using a single sequence of input-output data generated by $\Sigma$, of which the input sequence is PE. The notion of input-output data is formalized below.
\begin{definition}
    A pair of real-valued length-$N$ sequences $(u_{[0,N-1]},y_{[0,N-1]})$ is called an input-output trajectory of $\Sigma$ in~\eqref{eq:td-system}, if there exists a state sequence $x_{[0,N-1]}$ satisfying~\eqref{eq:td-system-state} for $k\in\mathbb{Z}_{[0,N-2]}$ and~\eqref{eq:td-system-output} for $k\in\mathbb{Z}_{[0,N-1]}$. In that case, the triplet $(u_{[0,N-1]},x_{[0,N-1]},y_{[0,N-1]})$ is called an input-state-output trajectory of $\Sigma$.
\end{definition}

Next, we recall WFL, which was originally published in~\cite{Willems2005}.
\begin{lemma}\label{lem:td-WFL}
    Let $(\hat{u}_{[0,N-1]},\hat{x}_{[0,N-1]},\hat{y}_{[0,N-1]})$ be an input-state-output trajectory\footnote{We use the accent \^{} to denote data that was collected previously/off-line. For example, $\hat{y}_{[0,M-1]}$ denotes the sequence of outputs collected off-line.} of $\Sigma$ in~\eqref{eq:td-system} satisfying Assumption~\ref{asm:ctrb}. Suppose that $\hat{u}_{[0,N-1]}$ is PE of order $L+n_x$ with $L\in\mathbb{Z}_{\geqslant 1}$. Then, the following statements hold:
    \begin{enumerate}[label=(\roman*)]
        \item\label{item:td-wfl-1} The matrix \[\begin{bmatrix}H_1(x_{[0,N-L]})\\ H_L(u_{[0,N-1]})\end{bmatrix}\] has full row rank;
        \item\label{item:td-wfl-2} The pair of real-valued sequences $(u_{[0,L-1]},y_{[0,L-1]})$ of length $L$ is an input-output trajectory of $\Sigma$, if and only if there exists $g\in\mathbb{R}^{N-L+1}$ such that
        \begin{equation*}
            \def\arraystretch{\arraystretchval}
\left[\begin{array}{@{}c@{}}
                u_{[0,L-1]}\\\hdashline[2pt/2pt]
                y_{[0,L-1]}
            \end{array}\right] = \def\arraystretch{\arraystretchval}
\left[\begin{array}{@{}c@{}}
                H_L(\hat{u}_{[0,N-1]})\\\hdashline[2pt/2pt]
                H_L(\hat{y}_{[0,N-1]})
            \end{array}\right]g.
        \end{equation*}
    \end{enumerate}~
\end{lemma}
Lemma~\ref{lem:td-WFL} is a combination of~\cite[Corollary~2.(\romannumeral 3)]{Willems2005}, which is Lemma~\ref{lem:td-WFL}.\ref{item:td-wfl-1}, and~\cite[Theorem 1]{Willems2005}, which is Lemma~\ref{lem:td-WFL}.\ref{item:fd-wfl-2}. Lemma~\ref{lem:td-WFL}.\ref{item:td-wfl-1} means that a powerful rank condition on a particular matrix containing both state and input data can be achieved by exciting $\Sigma$ with an input of sufficient order of PE, which has proved useful in, e.g., subspace identification~\cite{Verhaegen1992} and data-driven state-feedback synthesis~\cite{DePersis2020,vanWaarde2020,DePersis2021}. Lemma~\ref{lem:td-WFL}.\ref{item:td-wfl-2}, on the other hand, expresses all possible input-output trajectories of $\Sigma$ of length $L$ as a linear combination of time-shifted windows of a single input-output trajectory of length $N$ collected off-line using an input sequence of sufficiently-high PE order. This result has found numerous applications in system identification, see, e.g.,~\cite{Verhaegen1992,Markovsky2005,Markovsky2006}, it has also proved remarkably powerful in the context of data-driven \emph{predictive} control over the recent years, see, e.g.,~\cite{Coulson2019,Berberich2021,Verhoek2021b,Schmitz2022}. Rather than using a single persistently-exciting input-output trajectory, Lemma~\ref{lem:td-WFL} has also been extended to allow multiple input-output trajectories that are collectively but not necessarily individually PE~\cite{vanWaarde2020}. In the next section, we summarize how Lemma~\ref{lem:td-WFL} is used for predictive control.

\subsection{Data-driven predictive control}\label{sec:DeePC}
In DeePC~\cite{Coulson2019}, Lemma~\ref{lem:td-WFL} is used to obtain a prediction model directly in terms of (time-domain) data. We complete this preliminaries section by briefly recalling the working principles of DeePC. Let $v_{i,k}\in\mathbb{R}^{n_v}$, $i\in\mathbb{Z}_{[m,n]}$ with $m,n\in\mathbb{Z}$ and $n\geqslant m$, denote the prediction of $v_{k+i}$ made at time $k\in\mathbb{Z}$. We introduce the notation $v_{[m,n],k}=\{v_{i,k}\}_{i\in\mathbb{Z}_{[m,n]}}$ to refer to the sequence of such predictions made at time $k\in\mathbb{Z}$ or to refer to the vectorized predicted sequence, i.e., \[v_{[m,n],k}=\begin{bmatrix} v_{m,k}^\top & v_{m+1,k}^\top & \hdots & v_{n,k}^\top\end{bmatrix}^\top.\]

Since DeePC is not based on a particular realization of $\Sigma$, its state is fundamentally unknown. To ensure that the initial state of the prediction model is consistent with the internal state of $\Sigma$, an initial input-output trajectory $(u_{[k-\bar{T},k-1]},y_{[k-\bar{T},k-1]})$, $\bar{T}\in\mathbb{Z}_{\geqslant \ell_{\Sigma}}$, is prepended. Since we do not assume observability of $\Sigma$, only the observable part of the internal state is uniquely determined (up to a similarity transformation) by the initial input-output trajectory of length $\bar{T}\geqslant \ell_{\Sigma}$. Since the directions in the unobservable space and any similarity transformation do not affect the future output, the predicted output is uniquely determined and consistent with the initial trajectory.

Let $T\in\mathbb{Z}_{\geqslant 1}$ be the prediction horizon and let $(\hat{u}_{[0,N-1]},\hat{y}_{[0,N-1]})$ be an input-output trajectory of $\Sigma$, with $\hat{u}_{[0,N-1]}$ being PE of order $\bar{T}+T+n_x\geqslant \ell_{\Sigma}+T+n_x$. Then, at every $k\in\mathbb{Z}_{\geqslant 0}$, given the initial input-output trajectory $(u_{[k-\bar{T},k-1]},y_{[k-\bar{T},k-1]})$, DeePC~\cite{Coulson2019} solves the finite-horizon optimal control problem
\begin{equation}
	\begin{array}{@{}rl@{}} \underset{\mathclap{\substack{u_{[0,T-1],k},\\ y_{[0,T-1],k},\\g_k,\sigma_k}}}{\minimize} & \displaystyle\lambda_\sigma\|\sigma_k\|_1+\lambda_g\|g_k\|_1+\sum_{\mathclap{i\in\mathbb{Z}_{[0,T-1]}}} \ell(y_{i,k},u_{i,k}),\\
	\st & \def\arraystretch{\arraystretchval}
\left[\begin{array}{@{}c@{}}
        u_{[k-\bar{T},k-1]}\\[-1mm]
        u_{[0,T-1],k}\\\hdashline[2pt/2pt]
        y_{[k-\bar{T},k-1]}+\sigma_k\\[-1mm]
        y_{[0,T-1],k}
    \end{array}\right] = \def\arraystretch{\arraystretchval}
\left[\begin{array}{@{}c@{}}
        H_{\bar{T}+T}(\hat{u}_{[0,N-1]})\\\hdashline[2pt/2pt]
        H_{\bar{T}+T}(\hat{y}_{[0,N-1]})
    \end{array}\right]g_k,\\[9mm]
    & u_{i,k}\in\mathbb{U},~y_{i,k}\in\mathbb{Y},~\text{ for all }i\in\mathbb{Z}_{[0,T-1]},
    \end{array}\label{eq:DeePC}
\end{equation}
where $\mathbb{U}$ and $\mathbb{Y}$ denote the set of admissible inputs and outputs, respectively, while $\ell\colon\mathbb{R}^{n_y}\times\mathbb{R}^{n_u}\rightarrow\mathbb{R}_{\geqslant 0}$ denotes the stage cost. The decision variables in~\eqref{eq:DeePC} are $g_k\in\mathbb{R}^{N-T-\bar{T}+1}$, $\sigma_k\in\mathbb{R}^{\bar{T}n_y}$, and the predicted inputs and outputs $u_{i,k}$ and $y_{i,k}$, $i\in\mathbb{Z}_{[0,T-1]}$, respectively. Here, $\sigma_k$ is an auxiliary slack variable and $\lambda_\sigma,\lambda_g\in\mathbb{R}_{>0}$ are regularization parameters~\cite{Coulson2019}. These regularization parameters and slack variables are needed to deal with noise in the off-line data as well as the initial output trajectory $y_{[k-\bar{T},k-1]}$. In the noiseless/nominal case, however, we can use $\lambda_g=0$ and $\sigma_k=0$. To turn~\eqref{eq:DeePC} into a feedback policy, denote by $u_{[0,T-1],k}^*=\{u^*_{i,k}\}_{i\in\mathbb{Z}_{[0,T-1]}}$ the optimal control action computed at time $k$ by solving~\eqref{eq:DeePC}. DeePC implements the first element of $u^\star_{[0,T-1],k}$, i.e., $u_k=u^\star_{0,k}$, and solves~\eqref{eq:DeePC} at time $k+1$ using the updated initial input-output trajectory $(u_{[k-\bar{T}+1,k]},y_{[k-\bar{T}+1,k]})$, thereby creating a (receding-horizon) feedback policy.

\section{Problem setting}
\label{sec:problem-setting}
Throughout this paper, we consider the LTI system $\Sigma$ introduced in~\eqref{eq:td-system} satisfying Assumption~\ref{asm:ctrb}. We assume that $\Sigma$ itself, so $(A,B,C,D)$ in~\eqref{eq:td-system}, is unknown, and that, instead, we are given (sampled) frequency-domain data. To formalize this notion of frequency-domain data, we first recall that the spectrum $V~\colon~\mathbb{W}\rightarrow\mathbb{C}^{n_v}$, with $\mathbb{W}\coloneqq (-\pi,\pi]$, of the time-domain sequence $\{v_k\}_{k\in\mathbb{Z}}$, with $v_k\in\mathbb{R}^{n_v}$ for all $k\in\mathbb{Z}$, is obtained by taking the discrete-time Fourier transform (DTFT), i.e.,
\begin{equation*}
    V(\omega) = \sum_{k\in\mathbb{Z}} v_ke^{-j\omega k}. 
\end{equation*}
We also recall that, since $\{v_k\}_{k\in\mathbb{Z}}$ is real-valued, the spectrum $V$ is symmetric, i.e., $V(\omega)=V^*(-\omega)$ for all $\omega\in\mathbb{W}$.

Consider an equidistant frequency grid\footnote{This choice is made for notational convenience and our results readily extend to non-equidistant frequency grids (see also Remark~\ref{rem:non-equidistant}).} consisting of the $M\in\mathbb{Z}_{\geqslant 1}$ frequencies
\begin{equation}\label{eq:omega}
    \hat{\omega}_k = \frac{\pi k}{M}, \qquad k\in\mathbb{Z}_{[0,M-1]}.
\end{equation}
We denote the sequence of these frequencies as $\hat{\omega}_{[0,M-1]}\coloneqq \{\hat{\omega}_k\}_{k\in\mathbb{Z}_{[0,M-1]}}$, which satisfies $\hat{\omega}_{[0,M-1]}\subset\mathbb{W}$. Such a grid of (normalized) frequencies is commonly used in, e.g., frequency-domain identification literature~\cite{Pintelon2012}. We consider frequency-domain data that consists of sequences $\hat{V}_{[0,M-1]}=\{\hat{V}_k\}_{k\in\mathbb{Z}_{[0,M-1]}}$, with $\hat{V}_k=\hat{V}(\hat{\omega}_k)$ for all $k\in\mathbb{Z}_{[0,M-1]}$, of samples of such spectra at the frequencies in $\hat{\omega}_{[0,M-1]}$. Since $\hat{\omega}_0=0$, the frequency-domain data encountered throughout the paper satisfies $V_0\in\mathbb{R}^{n_v}$. Moreover, we exploit the symmetry of $\hat{V}$ to only store frequencies in $[0,\pi)$ in $\hat{\omega}_{[0,M-1]}$. In the remainder of this paper, we use lower case, e.g., $v_{[0,N-1]}$, to refer to real-valued time-domain sequences, while the upper case, e.g., $V_{[0,M-1]}$ is used to refer to complex-valued frequency-domain sequences. 

In particular, our frequency-domain data consists, depending on the precise application, of pairs of sequences $(\hat{U}_{[0,M-1]},\hat{Y}_{[0,M-1]})$ or $(\hat{U}_{[0,M-1]},\hat{X}_{[0,M-1]})$ containing samples of the input and, respectively, the corresponding output or state spectrum of $\Sigma$. To formalize this, we recall that the transfer function $H~\colon~\mathbb{C}\rightarrow\mathbb{C}^{n_y\times n_u}$ of $\Sigma$ is given by
\begin{equation*}
    H(z) = C(zI-A)^{-1}B+D,\quad z\in\mathbb{C}.
\end{equation*}
\begin{definition}\label{def:input-output-spectra}
    A pair of complex-valued length-$M$ sequences $(\hat{U}_{[0,M-1]},\hat{Y}_{[0,M-1]})$ is called an input-output spectrum of $\Sigma$, if $\hat{U}_0\in\mathbb{R}^{n_u}$, $\hat{Y}_0\in\mathbb{R}^{n_y}$, and there exists a state spectrum $\hat{X}_{[0,M-1]}$ with $\hat{X}_0\in\mathbb{R}^{n_x}$ satisfying
    \begin{subequations}
        \begin{align}
            e^{j\hat{\omega}_k}\hat{X}_k &= A\hat{X}_k + B\hat{U}_k,\label{eq:fd-system-state}\\
            \hat{Y}_k &= C\hat{X}_k + D\hat{U}_k\label{eq:fd-system-output}
        \end{align}
        \label{eq:fd-system}%
    \end{subequations}
    for all $k\in\mathbb{Z}_{[0,M-1]}$. In that case, the triplet $(\hat{U}_{[0,M-1]},\hat{X}_{[0,M-1]},\hat{Y}_{[0,M-1]})$ is called an input-state-output spectrum of $\Sigma$, and the pair $(\hat{U}_{[0,M-1]},\hat{X}_{[0,M-1]})$ is called an input-state spectrum of $\Sigma$.
\end{definition}
Definition~\ref{def:input-output-spectra} is also commonly used in subspace identification literature, see, e.g.,~\cite{McKelvey1996,Overschee1996,Cauberghe2006}. Note that Definition~\ref{def:input-output-spectra} does not include any transient phenomena. A practical and simple method to obtain data that adheres to Definition~\ref{def:input-output-spectra} is, as illustrated in Example~\ref{ex:noisy-simulation} and Section~\ref{sec:case-study-freepc}, by (pre-stabilizing the system, if needed, and) measuring sufficiently long for the transient to dampen out and, subsequently, discarding the initial part of the data where transient phenomena play a significant role. Recently, this restriction to steady-state data was relaxed in~\cite{Meijer2025-cdc-arxiv} to accommodate also non-steady-state data. Definition~\ref{def:input-output-spectra} can also be applied directly to FRF measurements $\{H(e^{j\hat{\omega}_k})\}_{k\in\mathbb{Z}_{[0,M-1]}}$, as detailed in Remark~\ref{rem:FRF-measurements} below. To see this, note that $\hat{X}_k$ can be eliminated from~\eqref{eq:fd-system} such that $(\hat{U}_{[0,M-1]},\hat{Y}_{[0,M-1]})$ is an input-output spectrum of $\Sigma$, if and only if 
\begin{equation*}
        \hat{Y}_k = H(e^{j\hat{\omega}_k})\hat{U}_k,
\end{equation*}
for all $k\in\mathbb{Z}_{[0,M-1]}$. In this case, a pre-processing step (e.g., using local polynomial methods~\cite{Pintelon2012}) can be used to remove transient phenomena from non-steady-state data, which does not necessarily require periodic excitation.
\begin{remark}[Incorporating FRF measurements]\label{rem:FRF-measurements}
    Let $\{H(e^{j\hat{\omega}_k})\}_{k\in\mathbb{Z}_{[0,M-1]}}$ be FRF measurements of a system $\Sigma$. If $\Sigma$ is single-input single-output (SISO), We can incorporate this frequency-domain data in Definition~\ref{def:input-output-spectra} by setting $\hat{Y}_k=H(e^{j\hat{\omega}_k})$ and $\hat{U}_k=1$ for all $k\in\mathbb{Z}_{[0,M-1]}$. The input-output spectrum $\{\hat{U}_{[0,M-1]},\hat{Y}_{[0,M-1]}\}$ then comprises our frequency-domain data. If $\Sigma$ is a MIMO system, however, we can only incorporate a single input direction $r_k\in\mathbb{C}^{n_u}$ for each $k\in\mathbb{Z}_{[0,M-1]}$ in Definition~\ref{def:input-output-spectra} by setting $\hat{Y}_k=H(e^{j\hat{\omega}_k})r_k$ and $\hat{U}_k = r_k$ for all $k\in\mathbb{Z}_{[0,M-1]}$. Clearly, this discards a lot of the available data. To make efficient use of all available data, we present a natural extension of our results (see Section~\ref{sec:collective-PE-WFL} below), where our data consists of a collection of $E$ input-output spectra $\{(\hat{U}^e_{[0,M-1]},\hat{X}^{e}_{[M-1]},\hat{Y}^{e}_{[0,M-1]})\}_{e\in\mathcal{E}}$ with $\mathcal{E}\coloneqq\mathbb{Z}_{[1,E]}$. This enables us to utilize all available data effectively also in the MIMO case by setting, e.g., $E = n_u$, $\hat{U}^e_k = \bm{e}_e$ and $\hat{Y}^e_k = H(e^{j\hat{\omega}_k})\hat{U}^e_k$ for all $e\in\mathcal{E}$ and $k\in\mathbb{Z}_{[0,M-1]}$. This way, we divide the data across $n_u$ data sets in each of which only one of the inputs is active.
\end{remark}
\begin{remark}\label{rem:non-equidistant}
    The results in the sequel also readily apply if the frequencies $\hat{\omega}_{[0,M-1]}$ lie on a non-equidistant grid which does not necessarily contain $\hat{\omega}_0=0$, as shown in~\cite{Meijer2024-ecc-arxiv}. In fact, to adhere to the conditions in Definition~\ref{def:input-output-spectra},~\cite{Meijer2024-ecc-arxiv} requires the frequency-domain data to be obtained from \emph{periodic} time-domain data, which implies the existence of the greatest common divisor of the excited frequencies. Hence, we can always define a grid as in~\eqref{eq:omega} (possibly containing many unexcited frequencies, for which $\hat{U}_k=0$, $\hat{X}_k=0$, and $\hat{Y}_k=0$) by choosing $M$ to be sufficiently large. 
\end{remark}

As indicated in the introduction, almost all existing variants of WFL use time-domain data and much less effort has been directed towards data-driven control based on WFL using frequency-domain data. In this paper, our goal is to fill the resulting gaps in the literature by formulating a version of WFL based on frequency-domain data and, subsequently, using it for data-driven (predictive) control.

\section{Willems' lemma in frequency domain}
\label{sec:frequency-domain-Willems}
In this section, we subsequently introduce the notion of persistence of excitation for frequency-domain data and use such data to formulate a frequency-domain version of Willems' fundamental lemma. Analogous to~\cite{vanWaarde2020}, we also present extensions of these results to exploit multiple frequency-domain data sets. 

\subsection{Persistence of excitation in frequency domain}
The original WFL requires the considered time-domain data to be sufficiently ``rich'', which is formalized in the notion of persistence of excitation (see Definition~\ref{def:td-PE}). It should come as no surprise that we also require that the frequency-domain data is sufficiently ``rich'', which we will now formalize by defining PE for the complex-valued sequence $V_{[0,M-1]}$ introduced in the previous section. To this end, let, for $L\in\mathbb{Z}_{\geqslant 1}$, $F_L\colon \mathbb{C}^{n_v(n-m+1)}\rightarrow\mathbb{C}^{n_vL\times (n-m+1)}$, where $m,n\in\mathbb{Z}_{\geqslant 0}$ with $n\geqslant m$, be the matrix-valued function of a sequence of length $n-m+1$ given by
\begin{equation}\label{eq:F_L}
    F_L(V_{[m,n]}) = \begin{bmatrix}
        W_L(e^{j\hat{\omega}_m})\otimes V_m & \hdots & W_L(e^{j\hat{\omega}_{n}})\otimes V_{n}
    \end{bmatrix},
\end{equation}
where $W_L(z)\coloneqq \begin{bmatrix} 1 & z & \hdots & z^{L-1}\end{bmatrix}^\top$, $z\in\mathbb{C}$. 

Using $F_L$ and the symmetry of the underlying spectrum, we define PE for the complex-valued sequence $V_{[0,M-1]}$ below.
\begin{definition}\label{def:fd-PE}
    The complex-valued sequence $V_{[0,M-1]}\in\mathbb{C}^{n_vM}$ of length $M$ with $V_0\in\mathbb{R}^{n_v}$ is said to be persistently exciting of order $L\in\mathbb{Z}_{[1,2M-1]}$, if the matrix 
    \begin{equation*}
        \Phi_L(V_{[0,M-1]})\coloneqq \begin{bmatrix} F_L(V_{[0,M-1]}) & F_L^*(V_{[1,M-1]})\end{bmatrix}
    \end{equation*} 
    has full row rank.
\end{definition}
Definition~\ref{def:fd-PE} exploits the symmetry in the underlying spectrum by including also the complex conjugate $F_L^*(V_{[1,M-1]})$ in $\Phi_L(V_{[0,M-1]})$. This term corresponds to the spectral content at the negative frequencies. Since $\hat{\omega}_0=0$ corresponds to the DC-gain of the system, the complex conjugate of $V_{0}$ is \emph{not} included. Persistence of excitation of order $L$ combined with the particular complex-conjugate structure of $\Phi_L(V_{[0,M-1]})$ means that we can express any length-$L$ sequence of $n_v$-dimensional \emph{real-valued} vectors as a linear combination of the columns of $\Phi_L(V_{[0,M-1]})$. Mathematically, this means that \begin{multline}\label{eq:Xcal}\mathcal{X}\coloneqq\Big\{\xi\in\mathbb{R}^{n_vL}~\Big|~\exists G_0\in\mathbb{R},G_1\in\mathbb{C}^{M-1} \text{ s.t. }\\\xi=\Phi_L(V_{[0,M-1]})(G_0,G_1,G_1^*)\Big\}=\mathbb{R}^{n_vL}.\end{multline} 
To see this, suppose that $V_{[0,M-1]}$ is PE of order $L$. Note that we can express $(G_0,G_1,G_1^*) = T_{\Re}(G_0,\Re G_1,-\Im G_1)$, where \begin{equation}
    T_{\Re} \coloneqq \frac{1}{2}\begin{bmatrix}
        2 & 0 & 0\\
        0 & I_{M-1} & -jI_{M-1}\\
        0 & I_{M-1} & jI_{M-1}
    \end{bmatrix}.\label{eq:T_Re}
\end{equation}
Substitution in~\eqref{eq:Xcal}, combined with the fact that $(G_0,\Re G_1,-\Im G_1)$ can be chosen arbitrarily in $\mathbb{R}^{2M-1}$, shows that \begin{multline*}\mathcal{X} = \im \Phi_L(V_{[0,M-1]})T_\Re =\\\im \begin{bmatrix} 
            \Re F_L(V_{[0,M-1]}) & \Im F_L(V_{[1,M-1]})
        \end{bmatrix}.\end{multline*}
Since $T_{\Re}$ is non-singular and, by the assumed PE, $\rank \Phi_L(V_{[0,M-1]})=n_vL$, we find that $\mathcal{X}=\mathbb{R}^{n_vL}$. Interestingly, it also follows that the rank condition in Definition~\ref{def:fd-PE} is equivalent to the \emph{real-valued} matrix \begin{equation*}
    \Phi_L(V_{[0,M-1]})T_{\Re} = \begin{bmatrix} 
            \Re F_L(V_{[0,M-1]}) & \Im F_L(V_{[1,M-1]})
        \end{bmatrix}
\end{equation*} having full row rank.

The rank condition in Definition~\ref{def:fd-PE} implies that $2M-1\geqslant Ln_v$, which provides a condition on the minimum amount of frequency-domain data that is required for a PE of order $L$. Here, data at any non-zero frequency, i.e., any frequency in $\hat{\omega}_{[1,M-1]}$, yields at most $2$ orders of PE because the complex conjugate also contributes, whereas the data at $\hat{\omega}_0=0$ contributes at most $1$ order of PE. In Section~\ref{sec:collective-PE-WFL}, we extend Definition~\ref{def:fd-PE} to allow the combination of multiple frequency-domain data sets that are not PE when considered individually but nonetheless are collectively sufficiently ``rich''.

\begin{remark}\label{rem:td-fd-equiv-PE}
    There is a close connection between persistence of excitation in the time domain (Definition~\ref{def:td-PE}) and in the frequency domain (Definition~\ref{def:fd-PE}). In fact, it can be shown that the length-$2M$ frequency-domain sequence $V_{[0,2M-1]}$ is PE of order $L$ if and only if the time-domain sequence $v_{[2M-1]}$, which is obtained by applying the inverse DTFT, i.e., \[v_k = \frac{1}{2\pi}\int_{-\pi}^{\pi} V(\omega)e^{j\omega k}\,\mathrm{d}\omega, \quad k\in\mathbb{Z}_{[0,2M-1]},\] to the spectrum \[V(\omega)=\sum_{\mathclap{k\in\mathbb{Z}_{[0,M-1]}}}V_k\delta(\omega-\hat{\omega}_k),\quad V(-\omega)=V^*(\omega),\quad \omega\in[0,\pi],\] is PE of order $L$ in the sense of Definition~\ref{def:fd-PE}.
\end{remark}

\subsection{A frequency-domain Willems' fundamental lemma}\label{sec:FD-WFL}
Next, we introduce our frequency-domain version of WFL~\cite{Willems2005}, which provides counterparts to both statements in Lemma~\ref{lem:td-WFL} based on frequency-domain data. 
\begin{theorem}\label{thm:fd-WFL}
    Let $(\hat{U}_{[0,M-1]},\hat{X}_{[0,M-1]},\hat{Y}_{[0,M-1]})$ be an input-state-output spectrum of $\Sigma$ in~\eqref{eq:td-system} satisfying Assumption~\ref{asm:ctrb}. Suppose that $\hat{U}_{[0,M-1]}$ is PE of order $L+n_x$. Then, the following statements hold:
    \begin{enumerate}[label=(\roman*)]
    \item\label{item:fd-wfl-1} The matrix \[\begin{bmatrix}
        F_1(\hat{X}_{[0,M-1]}) & F^*_1(\hat{X}_{[1,M-1]})\\
        F_L(\hat{U}_{[0,M-1]}) & F^*_L(\hat{U}_{[1,M-1]})
    \end{bmatrix}\] has full row rank;
    \item\label{item:fd-wfl-2} The pair of real-valued (time-domain) sequences $(u_{[0,L-1]},y_{[0,L-1]})$ of length $L$ is an input-output trajectory of $\Sigma$, if and only if there exist a scalar $G_0\in\mathbb{R}$ and a vector $G_1\in\mathbb{C}^{M-1}$ such that
    \begin{equation*}
        \def\arraystretch{\arraystretchval}
\left[\begin{array}{@{}c@{}}
            u_{[0,L-1]}\\\hdashline[2pt/2pt]
            y_{[0,L-1]}
        \end{array}\right] = \def\arraystretch{\arraystretchval}
\left[\begin{array}{@{}c;{2pt/2pt}c@{}}
            F_L(\hat{U}_{[0,M-1]}) & F_L^*(\hat{U}_{[1,M-1]})\\\hdashline[2pt/2pt]
            F_L(\hat{Y}_{[0,M-1]}) & F_L^*(\hat{Y}_{[1,M-1]})
        \end{array}\right]\def\arraystretch{\arraystretchval}
\left[\begin{array}{@{}c@{}}
            G_0\\
            G_1\\\hdashline[2pt/2pt]
            G_1^*
        \end{array}\right].
    \end{equation*}
    \end{enumerate}
\end{theorem}\vspace*{\belowdisplayskip}
The proof of Theorem~\ref{thm:fd-WFL} can be found in the Appendix. Theorem~\ref{thm:fd-WFL}.\ref{item:fd-wfl-1} states that one can render a specific matrix, which contains samples of the state and input spectrum, full row rank by injecting an input spectrum whose order of PE is at least equal to $n_x+n_uL$. This result, which provides a frequency-domain counterpart of Lemma~\ref{lem:td-WFL}.\ref{item:fd-wfl-1}, is used in proving Theorem~\ref{thm:fd-WFL}.\ref{item:fd-wfl-2}. Moreover, its time-domain version has found numerous applications in, for instance, subspace identification, see, e.g.,~\cite{Verhaegen1992}, and state-feedback synthesis, see, e.g.,~\cite{DePersis2020,vanWaarde2020,DePersis2021}. In these applications, Theorem~\ref{thm:fd-WFL}.\ref{item:fd-wfl-1} can be used to perform these tasks directly based on frequency-domain data. In fact, in Section~\ref{sec:fd-lqr}, we showcase such an application to obtain an LQR directly in terms of the frequency-domain data. Theorem~\ref{thm:fd-WFL}.\ref{item:fd-wfl-2}, on the other hand, provides a counterpart to Lemma~\ref{lem:td-WFL}.\ref{item:fd-wfl-2} based on frequency-domain data, and states that we can characterize all input-output trajectories of $\Sigma$ of length $L$ using samples of its input-output spectrum collected off-line. Importantly, Theorem~\ref{thm:fd-WFL}.\ref{item:fd-wfl-2} can be applied \emph{without} measuring the state spectrum $\hat{X}_{[0,M-1]}$ as it only uses of the input-output spectrum. 

Some further remarks regarding Theorem~\ref{thm:fd-WFL} are in order. 
\begin{itemize}[leftmargin=*]
    \item We see that, analogous to the shift theorem of the DTFT~\cite{Proakis1996}, the time delays that are present in the Hankel matrices in Lemma~\ref{lem:td-WFL} translate to linear phase terms, i.e., the increasing powers of $e^{j\hat{\omega}_k}$ in $W_L(e^{j\hat{\omega}_k})$ in Theorem~\ref{thm:fd-WFL};
    \item Using conjugate symmetry and~\eqref{eq:T_Re}, Theorem~\ref{thm:fd-WFL}.\ref{item:fd-wfl-1} is equivalent to the real-valued matrix \begin{equation*}
        \begin{bmatrix}
            \Phi_1(\hat{X}_{[0,M-1]})\\
            \Phi_L(\hat{U}_{[0,M-1]})
        \end{bmatrix}T_{\Re} = \begin{bmatrix}
            \Re F_1(\hat{X}_{[0,M-1]}) & \Im F_1(\hat{X}_{[1,M-1]})\\
            \Re F_L(\hat{U}_{[0,M-1]}) & \Im F_L(\hat{U}_{[1,M-1]})
        \end{bmatrix}
    \end{equation*}
    having full row rank. Similarly, we obtain from Theorem~\ref{thm:fd-WFL}.\ref{item:fd-wfl-2} that the pair of real-valued length-$L$ sequences $(u_{[0,L-1]},y_{[0,L-1]})$ is an input-output trajectory of $\Sigma$, if and only if there exists $g\in\mathbb{R}^{2M-1}$ such that \begin{equation}\label{eq:rv-fd-wfl-2}\def\arraystretch{\arraystretchval}
\left[\begin{array}{@{}c@{}}
        u_{[0,L-1]}\\\hdashline[2pt/2pt]
        y_{[0,L-1]}
        \end{array}\right] = \def\arraystretch{\arraystretchval}
\left[\begin{array}{@{}c;{2pt/2pt}c@{}}
            \Re F_L(\hat{U}_{[0,M-1]}) & \Im F_L(\hat{U}_{[1,M-1]})\\\hdashline[2pt/2pt]
            \Re F_L(\hat{Y}_{[0,M-1]}) & \Im F_L(\hat{Y}_{[1,M-1]})
        \end{array}\right]g,
    \end{equation} where we performed the change of variables $g = T_{\Re}^{-1}G = (G_0,2\Re(G_1),-2\Im(G_1))$. The real-valued formulation in~\eqref{eq:rv-fd-wfl-2} may be preferable because it avoids complex-valued computations, reduces the number of decision variables by a factor of two, and removes the need to enforce the particular structure $(G_0,G_1,G_1^*)$ using equality constraints.
\end{itemize} 

\subsection{Extension to multiple data sets}\label{sec:collective-PE-WFL}
In this section, we provide generalizations of Definition~\ref{def:fd-PE} and Theorem~\ref{thm:fd-WFL}, which are our frequency-domain notions of PE and WFL, respectively, to allow for multiple data sets. A similar extension has been proposed for the time-domain WFL in~\cite{vanWaarde2020} to deal with, e.g., missing data samples and unstable systems. For frequency-domain data, using multiple data sets is natural, if not necessary, because it enables us to incorporate information regarding multiple input directions at the same frequency. As discussed in Remark~\ref{rem:FRF-measurements}, choosing the number of data sets $E\in\mathbb{Z}_{\geqslant 1}$ equal to the number of inputs, i.e., $E=n_u$, allows us to use complete FRF measurements of MIMO transfer functions for our data, i.e, $\{H(e^{j\hat{\omega}_k})\}_{k\in\mathbb{Z}_{[0,M-1]}}$, whereas Theorem~\ref{thm:fd-WFL} can only use one column of $H(e^{j\hat{\omega}_k})$ for one input direction per frequency. 

To enable the use of data containing multiple input directions at the same frequency, we will now introduce the notion of collective persistence of excitation (CPE)~\cite{vanWaarde2020} for frequency-domain data. In particular, we consider collections of $E\in\mathbb{Z}_{\geqslant 1}$ sequences $\{V_{[0,M-1]}^e\}_{e\in\mathcal{E}}$ with $V_{[0,M-1]}^e=\{V^e_{k}\}_{k\in\mathbb{Z}_{[0,M-1]}}$ and $V_k^e=V^e(\hat{\omega}_k)$ for all $k\in\mathbb{Z}_{[0,M-1]}$ and $e\in\mathcal{E}\coloneqq \mathbb{Z}_{[1,E]}$. For notational simplicity, we assume that each of the sequences corresponds to the same sequence of frequencies $\hat{\omega}_{[0,M-1]}$, which was defined in Section~\ref{sec:problem-setting}. If all excited frequencies admit a greatest common divisor, this is without loss of generality because we can always define a common frequency grid and set $V_k^e=0$ for the unexcited frequencies per experiment (see also Remark~\ref{rem:non-equidistant}). Let, for $L\in\mathbb{Z}_{\geqslant 1}$, $\mathcal{F}_L~\colon~\mathbb{C}^{n_v(n-m+1)E}\rightarrow\mathbb{C}^{n_vL\times (n-m+1)E}$, where $m,n\in\mathbb{Z}_{\geqslant 0}$ with $n\geqslant m$, be the matrix-valued function of $E$ sequences of length $n-m+1$ given by \begin{multline*}
    \mathcal{F}_L(\{V^e_{[m,n]}\}_{e\in\mathcal{E}})=\\
\begin{bmatrix}
        W_L(e^{j\hat{\omega}_m})\otimes \begin{bmatrix}\begin{smallmatrix} V^1_m & \hdots & V^E_m\end{smallmatrix}\end{bmatrix} & \hdots & W_L(e^{j\hat{\omega}_n})\otimes \begin{bmatrix}\begin{smallmatrix} V^1_n & \hdots & V^E_n\end{smallmatrix}\end{bmatrix}
    \end{bmatrix},
\end{multline*}
which is essentially $F_L(V_{[m,n]})$ but for multiple data sets that are sorted into one block per frequency.
\begin{definition}\label{def:fd-PE-multi}
    The collection of complex-valued sequences $\{V^e_{[0,M-1]}\}_{e\in\mathcal{E}}$ of length $M$ with $V^e_0\in\mathbb{R}^{n_v}$, $e\in\mathcal{E}$, is said to be collectively persistently exciting (CPE) of order $L\in\mathbb{Z}_{[1,E(2M-1)]}$, if the matrix \begin{multline*}\Psi_L(\{V_{[0,M-1]}^e\}_{e\in\mathcal{E}})\coloneqq\\\begin{bmatrix}
        \mathcal{F}_L(\{V^e_{[0,M-1]}\}_{e\in\mathcal{E}}) & \mathcal{F}_L^*(\{V^e_{[1,M-1]}\}_{e\in\mathcal{E}})\\
    \end{bmatrix}\end{multline*} has full row rank.
\end{definition}
Collective persistence of excitation holds if one of the sequences $V_{[0,M-1]}^e$, $e\in\mathcal{E}$, is persistently exciting (in the sense of Definition~\ref{def:fd-PE}), but it is not necessary. In other words, CPE is significantly more flexible than PE of the individual sequences. Importantly, if $n_v>1$, CPE allows more than one direction to be considered per frequency. The latter is particularly useful in the context of MIMO systems, as illustrated in, e.g., Remark~\ref{rem:FRF-measurements}. If the data is readily provided, then $E$ is the number of data sets. However, during experiment design $E$ can be chosen to achieve a trade-off between the number of experiments that need to be carried out to collect data, and the number of independent input directions that can be excited per frequency. Thus, smaller $E<n_u$ means that we need to excite a broader range of frequencies to achieve CPE.

Next, we present a similar generalization of Theorem~\ref{thm:fd-WFL}.
\begin{theorem}\label{thm:fd-WFL-multi}
    Let $\{(\hat{U}^e_{[0,M-1]},\hat{X}^e_{[0,M-1]},\hat{Y}^e_{[0,M-1]})\}_{e\in\mathcal{E}}$ be a collection of input-state-output spectra of $\Sigma$ in~\eqref{eq:td-system} satisfying Assumption~\ref{asm:ctrb}. Suppose that $\{\hat{U}^e_{[0,M-1]}\}_{e\in\mathcal{E}}$ is CPE of order $L+n_x$. Then, the following statements hold:
    \begin{enumerate}[label=(\roman*)]
        \item\label{item:fd-wfl-multi-1} The matrix \[\begin{bmatrix}
            \mathcal{F}_1(\{\hat{X}^e_{[0,M-1]}\}_{e\in\mathcal{E}}) & \mathcal{F}_1^*(\{\hat{X}^e_{[1,M-1]}\}_{e\in\mathcal{E}})\\
            \mathcal{F}_L(\{\hat{U}^e_{[0,M-1]}\}_{e\in\mathcal{E}}) & \mathcal{F}_L^*(\{\hat{U}^e_{[1,M-1]}\}_{e\in\mathcal{E}})
        \end{bmatrix}\]
        has full row rank;
        \item\label{item:fd-wfl-multi-2} The pair of real-valued time-domain sequences $(u_{[0,L-1]},y_{[0,L-1]})$ of length $L$ is an input-output trajectory of $\Sigma$, if and only if there exist vectors $G_0\in\mathbb{R}^E$ and $G_1\in\mathbb{C}^{E(M-1)}$ such that $G=(G_0,G_1,G_1^*)$ satisfies \begin{equation*}\def\arraystretch{\arraystretchval}
\left[\begin{array}{@{}c@{}}
            u_{[0,L-1]}\\\hdashline[2pt/2pt]
            y_{[0,L-1]}
        \end{array}\right] = 
        \def\arraystretch{\arraystretchval}
\left[\begin{array}{@{}c@{}}
            \Psi_L(\{\hat{U}^e_{[0,M-1]}\}_{e\in\mathcal{E}})\\\hdashline[2pt/2pt]
            \Psi_L(\{\hat{Y}^e_{[0,M-1]}\}_{e\in\mathcal{E}})
        \end{array}\right]\def\arraystretch{\arraystretchval}
G.
        \end{equation*}
    \end{enumerate}
\end{theorem}\vspace*{\belowdisplayskip}
A sketch of proof based on the proof of Theorem~\ref{thm:fd-WFL} can be found in the Appendix. Analogous to~\cite[Theorem~2]{vanWaarde2020}, which extends Lemma~\ref{lem:td-WFL} to allow for multiple time-domain data sets, Theorem~\ref{thm:fd-WFL-multi} extends Theorem~\ref{thm:fd-WFL} to carefully combine multiple frequency-domain data sets. In the special case $E=1$, Theorem~\ref{thm:fd-WFL-multi} reduces to Theorem~\ref{thm:fd-WFL}. The comments made in Section~\ref{sec:frequency-domain-Willems} regarding Theorem~\ref{thm:fd-WFL} also apply to Theorem~\ref{thm:fd-WFL-multi}.

\section{Applications}
\label{sec:applications}
In this section, we present several interesting applications of Theorem~\ref{thm:fd-WFL}, and, in the case of multiple data sets, Theorem~\ref{thm:fd-WFL-multi}. The code for all examples can be found at \url{https://github.com/TomasMeijer/FD-WFL}.

\subsection{Frequency-domain data-driven simulation}\label{sec:simulation}
It is well-known that the {\it steady-state} response of an asymptotically stable system~\eqref{eq:td-system} can be simulated using FRF measurements. In this section, we exploit our frequency-domain version of WFL to simulate not only the steady-state behavior of $\Sigma$ \emph{but also its transient response}, which is, to the best of the authors' knowledge, not possible using existing techniques. Moreover, $\Sigma$ is also not required to be asymptotically stable here, although the off-line data may have to be collected in closed-loop with a pre-stabilizing controller, as we will illustrate in Section~\ref{sec:case-study-freepc}. As discussed in Section~\ref{sec:collective-PE-WFL}, it is natural to consider multiple data sets, i.e., multiple input directions at the same frequency, when dealing with MIMO FRF measurements and, as such, we use Theorem~\ref{thm:fd-WFL-multi} as the basis for our FRF-based simulation algorithm.
\begin{proposition}\label{prop:simulation}
    Let $\{(\hat{U}^e_{[0,M-1]},\hat{Y}^e_{[0,M-1]})\}_{e\in\mathcal{E}}$ be a collection of input-output spectra of $\Sigma$ in~\eqref{eq:td-system} satisfying Assumption~\ref{asm:ctrb}. Let $L,L_0\in\mathbb{Z}_{\geqslant 1}$ and suppose that $\{\hat{U}^e_{[0,M-1]}\}_{e\in\mathcal{E}}$ is CPE of order $L+L_0+n_x$. For any length-$L_0$ past input-output trajectory $(u_{[-L_0,-1]},y_{[-L_0,-1]})$ of $\Sigma$ and a (future) input sequence $u_{[0,L-1]}$, there exist vectors $G_0\in\mathbb{R}^E$ and $G_1\in\mathbb{C}^{E(M-1)}$ such that $G=(G_0,G_1,G_1^*)$ satisfies \begin{equation}\def\arraystretch{\arraystretchval}\left[\begin{array}{@{}c@{}}
        u_{[-L_0,L-1]}\\\hdashline[2pt/2pt]
        y_{[-L_0,-1]}
    \end{array}\right] = 
    \def\arraystretch{\arraystretchval}\left[\begin{array}{@{}c@{}}
        \Psi_{L_0+L}(\{\hat{U}^e_{[0,M-1]}\}_{e\in\mathcal{E}})\\\hdashline[2pt/2pt]
        \Psi_{L_0}(\{\hat{Y}^e_{[0,M-1]}\}_{e\in\mathcal{E}})
    \end{array}\right]G,\label{eq:sim-G}\end{equation} and $(u_{[-L_0,L-1]},y_{[-L_0,L-1]})$ with \begin{equation}
        y_{[-L_0,L-1]}=\label{eq:sim-y}\\
        \def\arraystretch{\arraystretchval}\left[\begin{array}{@{}c@{}}
        \Psi_{L_0+L}(\{\hat{Y}^e_{[0,M-1]}\}_{e\in\mathcal{E}})
    \end{array}\right]G
    \end{equation} is an input-output trajectory of $\Sigma$. If, in addition, $L_0\geqslant \ell_{\Sigma}$, then $y_{[0,L-1]}$ is unique.
\end{proposition}
Proposition~\ref{prop:simulation} can be used to perform data-driven simulation of the unknown system $\Sigma$ directly based on frequency-domain data of the system. To do so, we determine $G$ satisfying~\eqref{eq:sim-G} for the past input-output trajectory $(u_{[-L_0,-1]},y_{[-L_0,-1]})$ and the future input sequence $u_{[0,L-1]}$, which always exists due to the PE condition. Subsequently, we use the obtained $G$ to compute the future output sequence $y_{[0,L-1]}$ from~\eqref{eq:sim-y}. Since we do not assume observability of $\Sigma$, only the observable part of the internal state is uniquely determined (up to a similarity transformation) by the initial input-output trajectory of length $L_0\geqslant \ell_{\Sigma}$. Since the directions in the unobservable space and any similarity transformation do not affect the future output, it is uniquely determined and consistent with the initial trajectory. We can transform Proposition~\ref{prop:fd-eval} as described in Section~\ref{sec:FD-WFL} to obtain an equivalent real-valued formulation. 
\begin{example}[FRF-based simulation using noise-free data]\label{ex:noise-free-simulation}
We consider the linearized four-tank system~\eqref{eq:td-system} of~\cite{Raff2006} with
\begin{align}
    A &= \begin{bmatrix}
        0.921 & 0 & 0.041 & 0\\
        0 & 0.918 & 0 & 0.033\\
        0 & 0 & 0.924 & 0\\ 
        0 & 0 & 0 & 0.937
    \end{bmatrix},\label{eq:4TS}\\ B &= \begin{bmatrix}
        0.017 & 0.001\\
        0.001 & 0.023\\
        0 & 0.061\\
        0.072 & 0
    \end{bmatrix},~
    C = \begin{bmatrix}
        1 & 0 & 0 & 0\\
        0 & 1 & 0 & 0
    \end{bmatrix},\text{ and }D=0,\nonumber
\end{align}
Since the system has $n_u=2$ inputs, we perform $E=n_u=2$ experiments to generate the off-line data in which we excite $M=50$ frequencies $\hat{\omega}_{[0,M-1]}$ in~\eqref{eq:omega}. Since we are dealing with noise-free data, we can, for the purposes of this example, generate the off-line data consisting of input-output spectra $\{(\hat{U}^e_{[0,M-1]},\hat{Y}^e_{[0,M-1]})\}_{e\in\mathcal{E}}$ by setting $\hat{U}^e_k=\bm{e}_e$ and computing $\hat{Y}^e_k = H(e^{j\hat{\omega}_k})\hat{U}^e_k$ for all $k\in\mathbb{Z}_{[0,M-1]}$ and $e\in\mathcal{E}$. In practice, we would obviously not be able to generate the data in this way because it requires a-priori knowledge of the transfer functions of $\Sigma$. Therefore, in Example~\ref{ex:noisy-simulation}, we also collect data by conducting (simulated) experiments including measurements noise.

The observability index of $\Sigma$ with~\eqref{eq:td-system} is given by $\ell_{\Sigma} = 2$. To perform the simulation, we use an initial trajectory of length $L_0=12> \ell_\Sigma$ using initial state $x_{-L_0}=0$ and input sequence $u_{[-L_0,-1]}$ with $u_{-L_0}=(1,0.5)\cdot10^{2}$ and $u_{[-L_0+1,-1]}=0$. Afterward, we perform a simulation of length $L=30$ using Proposition~\ref{prop:simulation} using input sequence $u_{[0,L-1]}$ with $u_{0}=(-1,-2)\cdot10^{2}$ and $u_{[1,L-1]}=0$. The results are depicted in Fig.~\ref{fig:noise-free-simulation} along with the true output response of the system $y^{\mathrm{true}}_{[-L_0,L-1]}$. We see that, despite the system being unstable, which causes the outputs to grow exponentially, the obtained results closely resemble the true solution, i.e., $\|y_{[0,L-1]}-y^{\mathrm{true}}_{[0,L-1]}\|_2=2.9662\cdot 10^{-14}$ and $\|y_{[0,L-1]}-y^{\mathrm{true}}_{[0,L-1]}\|_2/\|y^{\mathrm{true}}_{[0,L-1]}\|_2 = 2.0514\cdot 10^{-15}$.
\begin{figure}[!tb]
    \setlength\fwidth{0.38\textwidth}
    \begin{subfigure}{\fwidth}
%
%
\definecolor{mycolor1}{named}{blue}%
\definecolor{mycolor2}{named}{red}%
\definecolor{mycolor3}{named}{blue}%
\definecolor{mycolor4}{named}{red}%
\begin{tikzpicture}

\begin{axis}[%
width=\fwidth,
height=0.28\fwidth,
at={(-0.2\fwidth,0\fwidth)},
scale only axis,
xmin=-12,
xmax=29,
xlabel style={font=\color{white!15!black}},
xlabel={k},
ymin=-4,
ymax=2,
ylabel style={font=\color{white!15!black}},
ylabel={$y_k$},
axis background/.style={fill=white},
xmajorgrids,
ymajorgrids
]
\addplot [color=mycolor1, line width=2.0pt, forget plot]
  table[row sep=crcr]{%
0	1.32175740594804\\
1	-0.630244040690218\\
2	-1.03222109456486\\
3	-1.36810771986862\\
4	-1.64573446279854\\
5	-1.87211494182422\\
6	-2.05352545688628\\
7	-2.19557716400301\\
8	-2.30328148967351\\
9	-2.3811093995724\\
10	-2.43304508137296\\
11	-2.46263455165941\\
12	-2.47302965138289\\
13	-2.46702785280107\\
14	-2.44710826297253\\
15	-2.41546417433919\\
16	-2.37403248143314\\
17	-2.32452025402479\\
18	-2.26842873084621\\
19	-2.20707497415515\\
20	-2.1416114036512\\
21	-2.07304340843053\\
22	-2.00224521760154\\
23	-1.92997419372684\\
24	-1.85688369826623\\
25	-1.78353466454287\\
26	-1.71040600132225\\
27	-1.63790393877491\\
28	-1.56637041829047\\
29	-1.49609061823271\\
};
\addplot [color=mycolor2, line width=2.0pt, forget plot]
  table[row sep=crcr]{%
0	1.72100921871637\\
1	-3.00397350663131\\
2	-2.88642447042746\\
3	-2.77040151733791\\
4	-2.65629062363212\\
5	-2.5444139152751\\
6	-2.43503693226817\\
7	-2.32837516951093\\
8	-2.22459996156139\\
9	-2.12384377263884\\
10	-2.02620494770864\\
11	-1.93175197546386\\
12	-1.84052730943471\\
13	-1.75255078927455\\
14	-1.66782270045706\\
15	-1.58632650714072\\
16	-1.50803128978469\\
17	-1.43289391620939\\
18	-1.36086097215949\\
19	-1.29187047502568\\
20	-1.2258533921941\\
21	-1.16273498349911\\
22	-1.10243598544082\\
23	-1.04487365317623\\
24	-0.989962674789217\\
25	-0.937615970976011\\
26	-0.887745392037763\\
27	-0.840262322939494\\
28	-0.795078206165211\\
29	-0.752104991162893\\
};
\addplot [color=black, dashed, line width=1.2pt, forget plot]
  table[row sep=crcr]{%
-12	0\\
-11	1.75\\
-10	1.7368\\
-9	1.715139\\
-8	1.6864077078\\
-7	1.651832071335\\
-6	1.61249046664444\\
-5	1.56932921092463\\
-4	1.52317655707965\\
-3	1.47475531199825\\
-2	1.42469420785577\\
-1	1.37353814396213\\
0	1.32175740594804\\
1	-0.630244040690217\\
2	-1.03222109456486\\
3	-1.36810771986862\\
4	-1.64573446279853\\
5	-1.87211494182422\\
6	-2.05352545688628\\
7	-2.19557716400301\\
8	-2.3032814896735\\
9	-2.3811093995724\\
10	-2.43304508137296\\
11	-2.4626345516594\\
12	-2.47302965138288\\
13	-2.46702785280106\\
14	-2.44710826297252\\
15	-2.41546417433919\\
16	-2.37403248143313\\
17	-2.32452025402478\\
18	-2.2684287308462\\
19	-2.20707497415514\\
20	-2.14161140365119\\
21	-2.07304340843052\\
22	-2.00224521760154\\
23	-1.92997419372683\\
24	-1.85688369826622\\
25	-1.78353466454287\\
26	-1.71040600132225\\
27	-1.63790393877491\\
28	-1.56637041829047\\
29	-1.49609061823271\\
};
\addplot [color=black, dashed, line width=1.2pt, forget plot]
  table[row sep=crcr]{%
-12	0\\
-11	1.25\\
-10	1.3851\\
-9	1.494153\\
-8	1.5802378884\\
-7	1.646121673584\\
-6	1.69428880098485\\
-5	1.72696783034683\\
-4	1.74615570450545\\
-3	1.75363982109949\\
-2	1.75101810041792\\
-1	1.73971722591938\\
0	1.72100921871636\\
1	-3.00397350663131\\
2	-2.88642447042746\\
3	-2.77040151733791\\
4	-2.65629062363212\\
5	-2.5444139152751\\
6	-2.43503693226817\\
7	-2.32837516951092\\
8	-2.22459996156139\\
9	-2.12384377263884\\
10	-2.02620494770863\\
11	-1.93175197546385\\
12	-1.84052730943471\\
13	-1.75255078927454\\
14	-1.66782270045706\\
15	-1.58632650714072\\
16	-1.50803128978468\\
17	-1.43289391620939\\
18	-1.36086097215948\\
19	-1.29187047502567\\
20	-1.22585339219409\\
21	-1.16273498349911\\
22	-1.10243598544082\\
23	-1.04487365317623\\
24	-0.989962674789218\\
25	-0.937615970976013\\
26	-0.887745392037762\\
27	-0.840262322939496\\
28	-0.795078206165211\\
29	-0.752104991162892\\
};
\end{axis}
\end{tikzpicture}%
            \caption{Using noise-free data.}
            \label{fig:noise-free-simulation}
    \end{subfigure}\\
    \begin{subfigure}{\fwidth}
%
%
\definecolor{mycolor1}{named}{blue}%
\definecolor{mycolor2}{named}{red}%
\definecolor{mycolor3}{named}{blue}%
\definecolor{mycolor4}{named}{red}%
\begin{tikzpicture}

\begin{axis}[%
width=\fwidth,
height=0.28\fwidth,
at={(-0.2\fwidth,0\fwidth)},
scale only axis,
xmin=-12,
xmax=29,
xlabel style={font=\color{white!15!black}},
xlabel={k},
ymin=-4,
ymax=2,
ylabel style={font=\color{white!15!black}},
ylabel={$y_k$},
axis background/.style={fill=white},
xmajorgrids,
ymajorgrids
]
\addplot [color=mycolor1, line width=2.0pt, forget plot]
  table[row sep=crcr]{%
0	1.31440685349183\\
1	-0.62806315634975\\
2	-1.02085581318057\\
3	-1.34332324994467\\
4	-1.6297844960114\\
5	-1.86235304574968\\
6	-2.02147635921999\\
7	-2.18600519653332\\
8	-2.27832435799221\\
9	-2.35012651924948\\
10	-2.42877360983897\\
11	-2.46468225124246\\
12	-2.46108230465535\\
13	-2.45148284660791\\
14	-2.45079959868184\\
15	-2.40555779849318\\
16	-2.36632637806154\\
17	-2.33892971598214\\
18	-2.26969305144026\\
19	-2.1982221301579\\
20	-2.13971069204883\\
21	-2.10893274286114\\
22	-1.99600228650301\\
23	-1.91674131598959\\
24	-1.85513695739525\\
25	-1.77376474775938\\
26	-1.70220655752939\\
27	-1.60877063032227\\
28	-1.58922030418824\\
29	-1.48325544696967\\
};
\addplot [color=mycolor2, line width=2.0pt, forget plot]
  table[row sep=crcr]{%
0	1.72491643574021\\
1	-3.02222735934059\\
2	-2.88501825604268\\
3	-2.75922087048918\\
4	-2.66340646096229\\
5	-2.53528816905256\\
6	-2.44048685478111\\
7	-2.35594098426573\\
8	-2.22550466514145\\
9	-2.12233743490094\\
10	-2.05431629390783\\
11	-1.92988231860137\\
12	-1.81427385145461\\
13	-1.75370939726318\\
14	-1.68969689784598\\
15	-1.60232631844882\\
16	-1.51713503485895\\
17	-1.46446012481811\\
18	-1.35337885997942\\
19	-1.299249457259\\
20	-1.24536084241915\\
21	-1.16032592087217\\
22	-1.11103986445976\\
23	-1.02772104277665\\
24	-0.968804148859504\\
25	-0.911031800284383\\
26	-0.922873296382833\\
27	-0.858200336455601\\
28	-0.782321780570419\\
29	-0.746693866138273\\
};
\addplot [color=black, dashed, line width=1.2pt, forget plot]
  table[row sep=crcr]{%
-12	0\\
-11	1.75\\
-10	1.7368\\
-9	1.715139\\
-8	1.6864077078\\
-7	1.651832071335\\
-6	1.61249046664444\\
-5	1.56932921092463\\
-4	1.52317655707965\\
-3	1.47475531199825\\
-2	1.42469420785577\\
-1	1.37353814396213\\
0	1.32175740594804\\
1	-0.630244040690217\\
2	-1.03222109456486\\
3	-1.36810771986862\\
4	-1.64573446279853\\
5	-1.87211494182422\\
6	-2.05352545688628\\
7	-2.19557716400301\\
8	-2.3032814896735\\
9	-2.3811093995724\\
10	-2.43304508137296\\
11	-2.4626345516594\\
12	-2.47302965138288\\
13	-2.46702785280106\\
14	-2.44710826297252\\
15	-2.41546417433919\\
16	-2.37403248143313\\
17	-2.32452025402478\\
18	-2.2684287308462\\
19	-2.20707497415514\\
20	-2.14161140365119\\
21	-2.07304340843052\\
22	-2.00224521760154\\
23	-1.92997419372683\\
24	-1.85688369826622\\
25	-1.78353466454287\\
26	-1.71040600132225\\
27	-1.63790393877491\\
28	-1.56637041829047\\
29	-1.49609061823271\\
};
\addplot [color=black, dashed, line width=1.2pt, forget plot]
  table[row sep=crcr]{%
-12	0\\
-11	1.25\\
-10	1.3851\\
-9	1.494153\\
-8	1.5802378884\\
-7	1.646121673584\\
-6	1.69428880098485\\
-5	1.72696783034683\\
-4	1.74615570450545\\
-3	1.75363982109949\\
-2	1.75101810041792\\
-1	1.73971722591938\\
0	1.72100921871636\\
1	-3.00397350663131\\
2	-2.88642447042746\\
3	-2.77040151733791\\
4	-2.65629062363212\\
5	-2.5444139152751\\
6	-2.43503693226817\\
7	-2.32837516951092\\
8	-2.22459996156139\\
9	-2.12384377263884\\
10	-2.02620494770863\\
11	-1.93175197546385\\
12	-1.84052730943471\\
13	-1.75255078927454\\
14	-1.66782270045706\\
15	-1.58632650714072\\
16	-1.50803128978468\\
17	-1.43289391620939\\
18	-1.36086097215948\\
19	-1.29187047502567\\
20	-1.22585339219409\\
21	-1.16273498349911\\
22	-1.10243598544082\\
23	-1.04487365317623\\
24	-0.989962674789218\\
25	-0.937615970976013\\
26	-0.887745392037762\\
27	-0.840262322939496\\
28	-0.795078206165211\\
29	-0.752104991162892\\
};
\end{axis}
\end{tikzpicture}%
        \caption{Using noisy data.}
        \label{fig:noisy-simulation}
    \end{subfigure}
    \definecolor{mycolor1}{named}{blue}%
    \definecolor{mycolor2}{named}{red}%
    \definecolor{mycolor3}{named}{blue}%
    \definecolor{mycolor4}{named}{red}%
    \caption{Data-driven simulation of a linearized four-tank system using different data sets. In both cases, the simulated outputs $y_1$ (\kern0.5pt\protect\tikz[scale=0.6]{\protect\draw[color=mycolor1, line width=1.2pt] (0,0.1)--(17pt,0.1);\protect\draw[color=white] (0,0)--(0.55,0)}) and $y_2$ (\kern0.5pt\protect\tikz[scale=0.6]{\protect\draw[color=mycolor2, line width=1.2pt] (0,0.1)--(17pt,0.1);\protect\draw[color=white] (0,0)--(0.55,0)}) closely resemble the true response (\kern0.5pt\protect\tikz[scale=0.6]{\protect\draw[color=black, dashed, line width=1.2pt] (0,0.1)--(17pt,0.1);\protect\draw[color=white] (0,0)--(0.55,0)}).}
\end{figure}
\end{example}

\begin{example}[Data-driven simulation using noisy data]\label{ex:noisy-simulation}
We consider again the linearized four-tank system~\eqref{eq:4TS}. However, we now consider a more realistic setting with noisy data. To generate this data, we replace~\eqref{eq:td-system-output} by \[y_k = Cx_k+Du_k + n_k,\] to incorporate measurements noise. We conduct $E=n_u=2$ experiments to collect our off-line data, in each of which we excite only one direction of the input $u_k$ using a multi-sine. For experiment $e\in\mathcal{E}$, we use input $u^e_k$ given by
\begin{equation*}
    u_k^e = \begin{cases}
        (10\sum_{m\in\mathbb{Z}_{[0,M-1]}}\cos(\hat{\omega}_mk+\hat{\phi}^e_m),0),\text{ if }e=1,\\
        (0,10\sum_{m\in\mathbb{Z}_{[0,M-1]}}\cos(\hat{\omega}_mk+\hat{\phi}^e_m)),\text{ if }e=2,
    \end{cases}%
\end{equation*}
with $M=50$ frequencies $\hat{\omega}_{[0,M-1]}$ in~\eqref{eq:omega} and random phase shifts $\{\hat{\phi}^e_k\}_{k\in\mathbb{Z}_{[0,M-1]},e\in\mathcal{E}}$. Here, $\{n_k\}_{k\in\mathbb{Z}}$ is zero-mean white Gaussian noise with $\mathbb{E}[n_kn_k^\top]=\sqrt{0.1}I$ for all $k\in\mathbb{Z}$. We measure $u_k$ and $y_k$, of which we take the DFT to obtain frequency-domain input-output data. We collect our data by measuring $p+p_0$ periods of our multi-sine excitation, after which we discard the first $p_0=20$ periods to eliminate most of the transient effects. This elementary method, which approximates measurements of the steady-state behavior~\cite{Pintelon2012}, ensures that we (approximately) adhere to the conditions in Definition~\ref{def:input-output-spectra}. For $p=50$ periods, we implement the data-driven simulation based on Proposition~\ref{prop:simulation}, using the same initial trajectory of length $L_0=12>\ell_{\Sigma}$ and input sequence as in Example~\ref{ex:noise-free-simulation}, to simulate $L=30$ steps into the future. The results are shown in Fig.~\ref{fig:noisy-simulation} along with the true output response of the system $y^{\mathrm{true}}_{[-L_0,L-1]}$. To show the effect of measuring more periods, we also carried out this data-driven simulation for a varying numbers of periods $p$, and summarized the resulting absolute and relative errors in the simulated future outputs in Table~\ref{tab:sim}. For the sake of comparison, we also implemented a data-driven simulation using~\cite[Algorithm 1]{Markovsky2008} directly based on the underlying time-domain data and summarize the results in Table~\ref{tab:sim-td}. In both cases, it can be seen that the errors decrease as we use more periods, and that for higher values of $p$ the frequency-domain method becomes slightly more accurate. We emphasize that, when using frequency-domain data, this improvement does \emph{not} come at the cost of increased computational complexity because the dimensions of~\eqref{eq:sim-G} and~\eqref{eq:sim-y} do not scale with the number of periods $p$ but rather with the number of excited frequencies $M$ and with the number of data sets $E$. In the time-domain approach, on the other hand, increasing $p$ directly leads to increased computational cost because the dimensions of the involved Hankel matrices depend on the length of the time-domain data. These computational benefits are particularly interesting when using Proposition~\ref{prop:simulation} as the basis for a data-driven predictive control scheme, see Section~\ref{sec:data-driven-predictive-control} below, where computational cost is an important consideration due to the optimization that needs to be carried out in real time.

\begin{table}[!bt]
    \centering
    \caption{Absolute and relative error in simulated future output trajectory using frequency-domain data with $p$ periods.}
    \label{tab:sim}
    \begin{tabular}{@{}c|cc@{}}
        \hline
        \setlength{\tabcolsep}{3pt}
        $p$ & $\|y_{[0,L-1]}-y^\mathrm{true}_{[0,L-1]}\|_2$ & $\|y_{[0,L-1]}-y^\mathrm{true}_{[0,L-1]}\|_2/\|y^\mathrm{true}_{[0,L-1]}\|_2$\\\hline
        $5$ & $1.598\cdot 10^{-1}$ & $1.105\cdot 10^{-2}$\\
        $10$ & $1.402\cdot 10^{-1}$ & $9.697\cdot 10^{-3}$\\
        $50$ & $4.708\cdot 10^{-2}$ & $3.256\cdot 10^{-3}$\\
        $100$ & $3.833\cdot 10^{-2}$ & $2.651\cdot 10^{-3}$
    \end{tabular}
\end{table}
\begin{table}[!bt]
    \centering
    \caption{Absolute and relative error in simulated future output trajectory using time-domain data with $p$ periods.}
    \label{tab:sim-td}
    \begin{tabular}{@{}c|cc@{}}
        \hline
        \setlength{\tabcolsep}{3pt}
        $p$ & $\|y_{[0,L-1]}-y^\mathrm{true}_{[0,L-1]}\|_2$ & $\|y_{[0,L-1]}-y^\mathrm{true}_{[0,L-1]}\|_2/\|y^\mathrm{true}_{[0,L-1]}\|_2$\\\hline
        $5$ & $1.593\cdot 10^{-1}$ & $1.102\cdot 10^{-2}$\\
        $10$ & $1.405\cdot 10^{-1}$ & $9.717\cdot 10^{-3}$\\
        $50$ & $4.724\cdot 10^{-2}$ & $3.267\cdot 10^{-3}$\\
        $100$ & $3.942\cdot 10^{-2}$ & $2.726\cdot 10^{-3}$
    \end{tabular}
\end{table}
\end{example}

\subsection{Data-driven frequency response evaluation}
All the results presented so far mix time and frequency domain in the sense that they use frequency-domain off-line data to characterize the behavior of $\Sigma$ in the time domain. For example, both of the frequency-domain versions of WFL, i.e., Theorem~\ref{thm:fd-WFL} and~\ref{thm:fd-WFL-multi}, feature a time-domain response of $\Sigma$ on the left-hand side of Theorem~\ref{thm:fd-WFL}.\ref{item:fd-wfl-2} and Theorem~\ref{thm:fd-WFL-multi}.\ref{item:fd-wfl-multi-2}, respectively, and frequency-domain data on the respective right-hand sides. Then, we utilize Theorem~\ref{thm:fd-WFL-multi} to characterize the frequency-response of the system $\Sigma$ at any complex frequency that does not coincide with an eigenvalue of $A$.
\begin{proposition}\label{prop:fd-eval}
    Let $\{(\hat{U}^e_{[0,M-1]},\hat{Y}^e_{[0,M-1]})\}_{e\in\mathcal{E}}$ be a collection of input-output spectra of $\Sigma$ in~\eqref{eq:td-system} satisfying Assumption~\ref{asm:ctrb}. Let $L_0\in\mathbb{Z}_{\geqslant \ell_{\Sigma}}$, $L=L_0+1$, and suppose $\{\hat{U}^e_{[0,M-1]}\}_{e\in\mathcal{E}}$ is CPE of order $L_0 + 1 + n_x$. Then, for any sample $U_z\in\mathbb{C}^{n_u}$ of the input spectrum at $z\in\mathbb{C}$ which is not an eigenvalue of $A$, the system of linear equations
    \begin{equation}\label{eq:fd-eval}
        \def\arraystretch{\arraystretchval}\left[\begin{array}{@{}c;{2pt/2pt}c@{}}
            0 & \Psi_{L}(\{\hat{U}^{e}_{[0,M-1]}\}_{e\in\mathcal{E}})\\\hdashline[2pt/2pt]
            -W_{L}(z)\otimes I_{n_y} & \Psi_{L}(\{\hat{Y}^{e}_{[0,M-1]}\}_{e\in\mathcal{E}})
        \end{array}\right]\left[\begin{array}{@{}c@{}}
            Y_z\\\hdashline[2pt/2pt]
            G
        \end{array}\right] 
        \def\arraystretch{\arraystretchval}= \left[\begin{array}{@{}c@{}}
            \mathcal{U}_{L,z}\\\hdashline[2pt/2pt]
            0         
        \end{array}\right],
    \end{equation}
    where $\mathcal{U}_{L,z}=W_{L}(z)\otimes U_z$, has a unique solution for $Y_z$, which is such that the complex-valued pair $(U_z,Y_z)$ is a sample of the input-output spectrum of $\Sigma$ at $z$ (i.e., $Y_z = H(z)U_z$). 
\end{proposition}

Proposition~\ref{prop:fd-eval} allows us to evaluate the frequency response of $\Sigma$ to an input ``direction'' $U_z$ at the complex frequency $z$. Note that $z$ does not necessarily have to correspond to any of the frequencies $e^{j\hat{\omega}_k}$, $k\in\mathbb{Z}_{[0,M-1]}$, present in our off-line data, In fact, $z$ does not even have to be on the unit circle. Interestingly, if we take $U_z$ to be the identity matrix and apply Proposition~\ref{prop:fd-eval} to each of its columns individually, we can compute $Y_z=H(z)$ and, thereby, we can compute the transfer function (matrix) of $\Sigma$ at $z$. In fact, in this case, Proposition~\ref{prop:fd-eval} recovers a version of~\cite[Theorem 2]{Markovsky2024} that is using frequency-domain data. Importantly, unlike in the other results throughout this paper, $G$ in Proposition~\ref{prop:fd-eval} is not constrained to have a complex-conjugate structure (i.e., $(G_0,G_1,G_1^*)$ in, e.g., Theorem~\ref{thm:fd-WFL}.\ref{item:fd-wfl-2}). This is because the right-hand side of~\eqref{eq:fd-eval} is complex-valued.

\subsection{Frequency-domain data-driven LQR}\label{sec:fd-lqr}
In this section, we use Theorem~\ref{thm:fd-WFL} to perform a frequency-domain data-based LQR using input-state data, i.e., frequency-domain data consisting of input-state spectra. To be precise, we aim to design the optimal control law 
\begin{equation}\label{eq:lqr}
    u_k = Kx_k,
\end{equation}
for the system $\Sigma$ in~\eqref{eq:td-system} such that $\lim_{k\rightarrow \infty}x_k=0$, and such that the cost function
\begin{equation}\label{eq:lqr-cost}
    J(\{x_k\}_{k\in\mathbb{Z}_{\geqslant 0}},\{u_k\}_{k\in\mathbb{Z}_{\geqslant 0}}) = \sum_{k\in\mathbb{Z}_{\geqslant 0}} x_k^\top Qx_k + u_k^\top Ru_k,
\end{equation}
is minimized. In~\eqref{eq:lqr-cost}, $Q\in\mathbb{R}^{n_x\times n_x}$ is symmetric positive semi-definite and $R\in\mathbb{R}^{n_u\times n_u}$ is symmetric positive definite. 

\begin{proposition}\label{prop:lqr}
    Let $\{(\hat{U}^e_{[0,M-1]},\hat{X}^e_{[0,M-1]})\}_{e\in\mathcal{E}}$ be a collection of input-state spectra of $\Sigma$ in~\eqref{eq:td-system} satisfying Assumption~\ref{asm:ctrb}. Suppose that $\{\hat{U}^e_{[0,M-1]}\}_{e\in\mathcal{E}}$ is CPE of order $n_x+1$ and that every eigenvalue of $A$ on the unit circle is $(Q,A)$-observable, i.e., $\rank\begin{bmatrix} Q & \lambda I-A\end{bmatrix}= n_x$ for all $\lambda\in\mathbb{C}$ with $|\lambda|=1$. Let \[\def\arraystretch{\arraystretchval}\Delta \coloneqq \left[\begin{array}{@{}c@{}}
        X_0\\
        X_1\\\hdashline[2pt/2pt]
        U
    \end{array}\right] \coloneqq \left[\begin{array}{@{}cc@{}}
            \mathcal{F}_2(\{\hat{X}^e_{[0,M-1]}\}_{e\in\mathcal{E}}) & \mathcal{F}_2^*(\{\hat{X}^e_{[1,M-1]}\}_{e\in\mathcal{E}})\\\hdashline[2pt/2pt]
            \mathcal{F}_1(\{\hat{U}^e_{[0,M-1]}\}_{e\in\mathcal{E}}) & \mathcal{F}_1^*(\{\hat{U}^e_{[1,M-1]}\}_{e\in\mathcal{E}})
        \end{array}\right],\] and let 
            $\Pi(P)\coloneqq \diag\{Q-P,P,R\}$.
    Then, there exists a right inverse $X_0^\dagger$ of $X_0$, which is such that $\Delta^\hop\Pi(P)\Delta X_0^\dagger=0$ and such that the optimal control law that minimizes the cost~\eqref{eq:lqr-cost} and achieves $\lim_{k\rightarrow\infty}x_k=0$ is given by~\eqref{eq:lqr} and
    \begin{equation}\label{eq:K-opt}
        K = UX_0^\dagger,
    \end{equation}
    where $P$ is the unique solution to
    \begin{align}\label{eq:lqr-opt}
        \begin{array}{@{}rl@{}}
        \underset{\mathclap{P\in\mathbb{R}^{n_x\times n_x}}}{\maximize} & \operatorname{trace} P\\
        \st & \def\arraystretch{\arraystretchval}\Delta^\hop\Pi(P)\Delta\succcurlyeq 0, \quad P=P^\top\succcurlyeq 0.
        \end{array}
    \end{align}~
\end{proposition}
Proposition~\ref{prop:lqr} presents what is essentially a frequency-domain version of~\cite[Theorem~29]{vanWaarde2020c}, but we assume that our off-line data is persistently exciting instead of using the informativity approach. 
The right inverse $X_0^\dagger$ in~\eqref{eq:K-opt} is not necessarily unique, but we can construct one as in~\eqref{eq:X_0dag} below.

\begin{example}[Frequency-domain data-driven LQR]\label{ex:fd-lqr}
Consider again the linearized four-tank system given by~\eqref{eq:td-system} with $A$ and $B$ as in~\eqref{eq:4TS}, but with full state measurements, i.e., $C=I$ and $D=0$. We again generate our off-line frequency-domain data for the $M=50$ frequencies $\hat{\omega}_{[0,M-1]}$ in~\eqref{eq:omega} as outlined in Example~\ref{ex:noise-free-simulation}. We use Proposition~\ref{prop:lqr} to compute the LQR based on the resulting input-state spectra $\{(\hat{U}^e_{[0,M-1]},\hat{X}^e_{[0,M-1]})\}_{e\in\mathcal{E}}$. In particular, we minimize the cost function~\eqref{eq:lqr-cost} with $Q=I$ and $R=I$ by, first, solving~\eqref{eq:lqr-opt} to obtain $P$ and, subsequently, computing $K$ as in~\eqref{eq:K-opt} with 
\begin{equation}\label{eq:X_0dag}
    X_0^\dagger=(\Delta^\hop\Pi(P)\Delta)_{\perp}(X_0(\Delta^\hop\Pi(P)\Delta)_\perp)^+.
\end{equation}
Observe from~\eqref{eq:X_0dag} that $X_0^\dagger$ is, indeed, a right inverse of $X_0$ while also satisfying $\Delta^\hop\Pi(P)\Delta X_0^\dagger=0$. We find that
\begin{align*}
    P &= \begin{bmatrix}
        6.5015  & -0.0996 &   1.4346 &   -0.3019\\
        -0.0996  &  6.2416 &  -0.3062  &  1.0925\\
        1.4346  & -0.3062 &   6.7112  & -0.2291\\
        -0.3019   & 1.0925 &  -0.2291 &   7.0328
    \end{bmatrix},\\
    K &= \begin{bmatrix}-0.0785 &  -0.0733  & -0.0094 &  -0.4562\\
   -0.0820  & -0.1113   &-0.3667  & -0.0128\end{bmatrix}.
\end{align*}
As expected based on Proposition~\ref{prop:lqr}, these results closely resemble the true solution $\bar{P}$ to the discrete algebraic Riccati equation (DARE) and the true LQR gain $\bar{K}$, i.e., $\|\bar{P}-P\|_2=2.428\cdot 10^{-7}$, $\|\bar{P}-P\|_2/\|\bar{P}\|_2=2.886\cdot 10^{-8}$, $\|\bar{K}-K\|_2=1.429\cdot 10^{-8}$, and $\|\bar{K}-K\|_2/\|\bar{K}\|_2=2.995\cdot 10^{-8}$.
\end{example}

\section{Frequency-domain data-driven predictive~control}
\label{sec:data-driven-predictive-control}
In Section~\ref{sec:applications}, we have already seen how Theorem~\ref{thm:fd-WFL} and Theorem~\ref{thm:fd-WFL-multi} can be used to simulate the unknown system $\Sigma$ in~\eqref{eq:td-system} based on frequency-domain data. Naturally, we can also use these results as a frequency-domain data-driven substitute for the prediction model in model predictive control (MPC). The resulting frequency-domain data-driven predictive control (FreePC) scheme is the counterpart of DeePC~\cite{Coulson2019}.

\subsection{FreePC algorithm}
Let $T\in\mathbb{Z}_{\geqslant 1}$ be the length of the prediction horizon. We recall from Section~\ref{sec:simulation} that we need to use a past input-output sequence to ensure that our simulation, which we refer to as our prediction in the context of FreePC, starts from the correct initial condition. Let $\bar{T}\in\mathbb{Z}_{\geqslant \ell_{\Sigma}}$ denote the length of this past input-output sequence. Suppose we are given frequency-domain data consisting of the collection of $E$ input-output spectra $\{(\hat{U}^e_{[0,M-1]},\hat{Y}^e_{[0,M-1]})\}_{e\in\mathcal{E}}$ with $\{\hat{U}^e_{[0,M-1]}\}_{e\in\mathcal{E}}$ being PE of order $\bar{T}+T+n_x$. Note that this implies that $E(2M-1)\geqslant \bar{T}+T+n_x\geqslant \ell_{\Sigma}+T+n_x$. 

At every time $k\in\mathbb{Z}$, given the past input-output sequence $(u_{[k-\bar{T},k-1]},y_{[k-\bar{T},k-1]})$, the FreePC solution is now defined by solving the finite-horizon optimal control problem
\begin{align}
	\hspace*{.4cm}\setlength{\arraycolsep}{2pt}\begin{array}{@{}rl@{}} \underset{\mathclap{\substack{u_{[0,T-1],k},\\ y_{[0,T-1],k},\\G_{0,k},G_{1,k},\sigma_k}}}{\minimize} & \displaystyle\lambda_\sigma\|\sigma_k\|_1+\lambda_G\|G_k\|_1+\sum_{\mathclap{i\in\mathbb{Z}_{[0,T-1]}}} \ell(y_{i,k},u_{i,k}),\\
	\st & \def\arraystretch{\arraystretchval}
\left[\begin{array}{@{}c@{}}
        u_{[k-\bar{T},k-1]}\\[-1mm]
        u_{[0,T-1],k}\\\hdashline[2pt/2pt]
        y_{[k-\bar{T},k-1]}+\sigma_k\\[-1mm]
        y_{[0,T-1],k}
    \end{array}\right] = \def\arraystretch{\arraystretchval}
\left[\begin{array}{@{}c@{}}
        \Psi_{\bar{T}+T}(\{\hat{U}^e_{[0,M-1]}\}_{e\in\mathcal{E}})\\\hdashline[2pt/2pt]
        \Psi_{\bar{T}+T}(\{\hat{Y}^e_{[0,M-1]}\}_{e\in\mathcal{E}})
    \end{array}\right]G_k,\\[7mm]
    &\def\arraystretch{\arraystretchval}G_k=(G_{0,k},G_{1,k},G^*_{1,k}),\\
    & G_{0,k}\in\mathbb{R}^E, G_{1,k}\in\mathbb{C}^{E(M-1)},\\
    & u_{i,k}\in\mathbb{U},~y_{i,k}\in\mathbb{Y},~\text{ for all }i\in\mathbb{Z}_{[0,T-1]}.
    \end{array}\raisetag{12pt}\label{eq:FreePC}
\end{align}
Let $u_{[0,T-1],k}=\{u^\star_{i,k}\}_{i\in\mathbb{Z}_{[0,T-1]}}$ denote the optimal control action computed at time $k\in\mathbb{Z}$ by solving~\eqref{eq:FreePC}. After solving~\eqref{eq:FreePC}, we implement the first element of the optimal control action $u^\star_{[0,T-1],k}$, i.e., $u_k = u^\star_{0,k}$, and, at time $k+1$, we solve~\eqref{eq:FreePC} again using the newly-obtained past input-output trajectory $(u_{[k-\bar{T}+1,k]},y_{[k-\bar{T}+1,k]})$. This way, the optimal control is computed iteratively using the receding-horizon principle and a feedback policy is created.

We recall that $u_{[0,T-1],k}$ and $y_{[0,T-1],k}$ denote the (vectorized) sequence of predictions of $u_{k+i}$ and $y_{k+i}$, $i\in\mathbb{Z}_{[0,T-1]}$, respectively, as introduced in Section~\ref{sec:DeePC}. Moreover, similar to DeePC, $\ell$ denotes the stage cost, $\mathbb{U}$ and $\mathbb{Y}$ denote, respectively, the sets of admissible inputs and outputs, $\sigma_k\in\mathbb{R}^{\bar{T}n_y}$ is an auxiliary slack variable, and $\lambda_\sigma,\lambda_G\in\mathbb{R}_{>0}$ are regularization parameters. The slack variable and the regularization parameters are needed to deal with noise in $y_{[k-\bar{T},k-1]}$ and in the data $\{(\hat{U}^e_{[0,M-1]},\hat{Y}^e_{[0,M-1]})\}_{e\in\mathcal{E}}$. Some further remarks regarding the optimization problem in~\eqref{eq:FreePC} are in order. 
\begin{itemize}[leftmargin=*]
\item The prediction model in~\eqref{eq:FreePC} involves the product between a complex-valued matrix and vector. Importantly, $G_{0,k}\in\mathbb{R}^E$ is real and $G_{1,k}\in\mathbb{C}^{E(M-1)}$ is complex, and the conjugate symmetry in~\eqref{eq:FreePC} ensures that the left-hand side of the prediction model is always real-valued. Depending on the solver that is being used, however, it is often beneficial to apply the change of variables involving $T_{\Re}$, as described in Section~\ref{sec:FD-WFL}, to obtain a real-valued formulation that is equivalent to~\eqref{eq:FreePC}.
\item In DeePC, the number of decision variables, through $g_k\in\mathbb{R}^{N-\bar{T}-T+1}$ in~\eqref{eq:DeePC}, depends on the length $N$ of the time-domain data. Interestingly, in FreePC, the number of decision variables, through $G_{0,k}$ and $G_{1,k}$ in~\eqref{eq:FreePC}, depends only on the number of frequencies $M$ and experiments $E$ and not on the length of the time-domain data used to compute $(\hat{U}_{[0,M-1]},\hat{Y}_{[0,M-1]})$. As we saw in Example~\ref{ex:noisy-simulation} and will see also in Section~\ref{sec:case-study-freepc}, we can exploit this to perform longer experiments and, thereby, reduce the effect of noise \emph{without} increasing the computational complexity of FreePC.
\end{itemize}

\subsection{An equivalence result}
Interestingly, when $\lambda_g=0$ in~\eqref{eq:DeePC} and $\lambda_G=0$ in~\eqref{eq:FreePC}, which is possible when the off-line data is noise-free, we can show that FreePC is equivalent to DeePC.
\begin{theorem}\label{thm:equivalence}
    Let $\{(\hat{U}^e_{[0,M-1]},\hat{Y}^e_{[0,M-1]})\}_{e\in\mathcal{E}}$ and $(\hat{u}_{[0,N-1]},\hat{y}_{[0,N-1]})$ be, respectively, a collection of input-output spectra and a trajectory of $\Sigma$ in~\eqref{eq:td-system} satisfying Assumption~\ref{asm:ctrb}. Let $\bar{T}\in\mathbb{Z}_{\geqslant \ell_{\Sigma}}$ and $T\in\mathbb{Z}_{\geqslant 1}$ and suppose that $\{\hat{U}^e_{[0,M-1]}\}_{e\in\mathcal{E}}$ and $\hat{u}_{[0,N-1]}$ are PE of order $\bar{T}+T+n_x$. Suppose that $\{(\hat{U}^e_{[0,M-1]},\hat{Y}^e_{[0,M-1]})\}_{e\in\mathcal{E}}$ and $(\hat{u}_{[0,N-1]},\hat{y}_{[0,N-1]})$ are noise-free and, accordingly, set $\lambda_g=0$ in~\eqref{eq:DeePC} and $\lambda_G=0$ in~\eqref{eq:FreePC}. Let $(u_{[k-\bar{T},k-1]},y_{[k-\bar{T},k-1]})$ be the past input-output trajectory available at time $k\in\mathbb{Z}$. Then, FreePC based on~\eqref{eq:FreePC} is equivalent to DeePC based on~\eqref{eq:DeePC} in the sense that, at any $k\in\mathbb{Z}$, the open-loop optimal control problems in~\eqref{eq:FreePC} and~\eqref{eq:DeePC} admit the same (set of) optimal control action(s) and corresponding output sequence(s).
\end{theorem}

It immediately follows from Theorem~\ref{thm:equivalence} that, if~\eqref{eq:DeePC} admits a unique optimal control action $u^{\star}_{0,k}$, then $u^\star_{0,k}$ is also the unique optimal control action obtained by solving~\eqref{eq:FreePC} and vice versa. The frequency-domain data $\{(\hat{U}^e_{[0,M-1]},\hat{Y}^e_{[0,M-1]})\}_{e\in\mathcal{E}}$ and time-domain data $(\hat{u}_{[0,N-1]},\hat{y}_{[0,N-1]})$ do \emph{not} need to be related through a Fourier transform and can be obtained from separate experiments. In fact, Theorem~\ref{thm:equivalence} considers FreePC based on {\it multiple} experiments and DeePC based on {\it a single} experiment, however, a similar equivalence holds for DeePC based on multiple experiments.

\begin{figure}[!tb]
    \centering
    \begin{tikzpicture}
        \node[minimum height=0.875cm, minimum width=1.375cm, fill=black!10, draw] (H) at (0,0) {$H(z)$};
        \node[minimum height=0.875cm, minimum width=1.375cm, fill=black!10, draw, left of=H, node distance=3cm] (C) {$-C(z)$};
        \node[circle, minimum size=4.5mm, draw, inner sep=0pt] at ($(H)!0.5!(C)$) (din)  {$+$};
        \node[circle, minimum size=4.5mm, draw, inner sep=0pt, right of=H, node distance=1.5cm] (nin)  {$+$};
        \draw[->] (din)--(H) node[pos=0.35, above, anchor=south] {$u_k$};
        \draw[->] (H)--(nin);
        \draw[<-] (din)--++(90:0.8cm) node[near end, left] {$d_k$};
        \draw[<-] (nin)--++(90:0.8cm) node[near end, left] {$n_k$};
        \draw[->,line cap=rect] (nin)--++(0:1cm) node[near end, above,anchor=south] {$y_k$};
        \draw[->] (nin)--++(0:0.5cm)--++(-90:0.8cm)-|($(C)+(180:1.5cm)+(-90:0.8cm)$)|-(C);
        \draw[->] (C)--(din);
    \end{tikzpicture}
    \caption{Closed-loop measurement setup.}
    \label{fig:closed-loop-meas}
\end{figure}
\subsection[Numerical case study]{Numerical case study}\label{sec:case-study-freepc}
The code for this numerical case study is available at \url{https://github.com/TomasMeijer/FD-WFL}. We consider the unstable batch reactor~\cite{Walsh2001,vanWaarde2020}, which we discretize using a sampling time of $T_s = 0.1$s to obtain a system of the form~\eqref{eq:td-system} with
\begin{align}
    A &= \begin{bmatrix}
         1.1782  &  0.0015 &   0.5116 &  -0.4033\\
        -0.0515  &  0.6619 &  -0.0110 &   0.0613\\
        0.0762   & 0.3351  &  0.5606  &  0.3824\\
        -0.0006  &  0.3353 &   0.0893 &   0.8494
    \end{bmatrix},\label{eq:batch-reactor}\\
    B &= \begin{bmatrix}
         0.0045  & -0.0876\\
    0.4672  &  0.0012\\
    0.2132  & -0.2353\\
    0.2131  & -0.0161
    \end{bmatrix},~ C^\top = \begin{bmatrix} 
        1 & 0\\
        0 & 1\\
        1 & 0\\
        -1 & 0
    \end{bmatrix},\text{ and }D=0.\nonumber
\end{align}
Since the system is unstable, we collect our data in a closed-loop measurement setup~\cite[Chapter~10]{Soderstrom1989}, as depicted in Fig,~\ref{fig:closed-loop-meas}, with the stabilizing controller \[
    C_{\mathrm{s}}(z) = \frac{1}{z-0.995}\begin{bmatrix} 0 & 1.09(z-0.900)\\ -2.717(z-0.840) & 0\end{bmatrix}.
\] 
Note that the tuning of such a controller may significantly impact the quality of the obtained data. In this closed-loop measurement setup, the input $u_k\in\mathbb{R}^{2}$ consists of the output of the controller on top of which we inject a signal $\{d_k\}_{k\in\mathbb{Z}}$ with $d_k\in\mathbb{R}^{2}$ for all $k\in\mathbb{Z}$. To incorporate noise in our case study, we replace the measurement model in~\eqref{eq:td-system-output} by \[y_k = Cx_k + Du_k + n_k,\] where $\{n_k\}_{k\in\mathbb{Z}}$ is zero-mean white Gaussian noise with $\mathbb{E}[n_kn_k^\top]=\sqrt{0.3}I$ for all $k\in\mathbb{Z}$.
We conduct $E=n_u=2$ experiments to collect our off-line data, in which we each excite only one element of the injected signal $d_k$ using a multi-sine. For experiment $e\in\mathcal{E}$, we inject the signal $d^e_k$ given by
\begin{equation*}
    d_k^e = \begin{cases}
        (10\sum_{m\in\mathbb{Z}_{[0,M-1]}}\cos(\hat{\omega}_mk+\hat{\phi}^e_m),0),\text{ if }e=1,\\
        (0,10\sum_{m\in\mathbb{Z}_{[0,M-1]}}\cos(\hat{\omega}_mk+\hat{\phi}^e_m)),\text{ if }e=2,
    \end{cases}
\end{equation*}
with the $M=40$ frequencies $\hat{\omega}_{[0,M-1]}$ in~\eqref{eq:omega} and random phase shifts $\{\hat{\phi}^e_k\}_{k\in\mathbb{Z}_{[0,M-1]},e\in\mathcal{E}}$. We measure $p+p_0$ periods of $u$ and $y$, which are both affected by the noise $n$ due to the closed-loop setting. By discarding the first $p_0=100$ periods, we get rid of most transient phenomena, after which we compute \emph{per period} the DTFT to obtain the frequency domain data $\{\hat{D}^{d}_\rho(\hat{\omega}_k)\}_{\rho\in\mathcal{P},k\in\mathbb{Z}_{[0,M-1]}}$, $\{\hat{U}^{d}_\rho(\hat{\omega}_k)\}_{\rho\in\mathcal{P},k\in\mathbb{Z}_{[0,M-1]}}$, and $\{\hat{Y}^{d}_\rho(\hat{\omega}_k)\}_{\rho\in\mathcal{P},k\in\mathbb{Z}_{[0,M-1]}}$, where $\mathcal{P}\coloneqq \mathbb{Z}_{[1,p]}$. We use this data to estimate the FRF of $\Sigma$ using the indirect method~\cite{vandenHof1996}, i.e., we first estimate the sensitivities, for all $k\in\mathbb{Z}_{[0,M-1]}$, \begin{align*}
    S_{yd}(\hat{\omega}_k) &\approx \frac{1}{p}\sum_{\rho\in\mathcal{P}} \hat{U}^d_\rho(\hat{\omega}_k)(\hat{D}^d_\rho(\hat{\omega}_k))^\hop,\\
    S_{ud}(\hat{\omega}_k) &\approx \frac{1}{p}\sum_{\rho\in\mathcal{P}} \hat{Y}^d_\rho(\hat{\omega}_k)(\hat{D}^d_\rho(\hat{\omega}_k))^\hop.
\end{align*}
Subsequently, we use the fact that \begin{align*}
    S_{yd}(\hat{\omega}_k) &= H(e^{j\hat{\omega}_k})(I+C_{\mathrm{s}}(e^{j\hat{\omega}_k})H(e^{j\hat{\omega}_k}))^{-1}D(\hat{\omega}_k)D^\hop(\hat{\omega}_k),\\
    S_{ud}(\hat{\omega}_k) &= (I+C_{\mathrm{s}}(e^{j\hat{\omega}_k})H(e^{j\hat{\omega}_k}))^{-1}D(\hat{\omega}_k)D^\hop(\hat{\omega}_k),
\end{align*}
to estimate the FRF as \[H(e^{j\hat{\omega}_k}) \approx S_{yd}(\hat{\omega}_k)S_{ud}^{-1}(\hat{\omega}_k) \text{ for all } k\in\mathbb{Z}_{[0,M-1]}.\] Figure~\ref{fig:frf} shows the resulting FRF measurement alongside their 99\% confidence intervals (approximated according to~\cite[Chapter 6]{vanBerkel2015}) using $p=10$ and $p=250$ periods. It is clear that increasing the number of periods significantly shrinks the confidence intervals and yields more accurate FRF measurements. Interestingly, in DeePC measuring more periods would result in a higher number of decision variables and, thereby, a higher computational complexity. To be precise, the number of decision variables in DeePC is $T(n_u+n_y)+\bar{T}n_y+N-T-\bar{T}+1$ with $N=2Mp$ being the duration of the time-domain data, which increases with the number of periods $p$. In FreePC, however, the number of decision variables, which is $T(n_u+n_y)+\bar{T}n_y+E(2M-1)$, remains the same regardless of the number of periods measured. Therefore, FreePC allows us to reduce the effect of noise by measuring more periods \emph{without} increasing the computational complexity of the resulting optimal control problem. This illustrates an important benefit of the frequency-domain WFL and FreePC over their respectively time-domain counterparts. As mentioned before, the ease of performing data collection in such a closed-loop setting is, particularly, for unstable systems, another important benefit of considering frequency-domain data. Doing so in the time domain is significantly less tractable because of the (de)convolutions involved, and in open loop the system trajectories would rapidly diverge causing numerical issues. For further details on time versus frequency domain identification, we refer the interested reader to~\cite{Schoukens2004}.
\begin{figure}[!tb]
    \setlength\fwidth{0.38\textwidth}
    \input{images/freepc-frf.tex}
    \caption{Estimated FRFs of the system (\kern0.5pt\protect\tikz[scale=0.6]{\protect\draw[dashed] (0,0.1)--(0.55,0.1);\protect\draw[color=white] (0,0)--(0.55,0)}) using $p=10$ (\kern0.5pt\protect\tikz[scale=0.3]{\protect\filldraw[color=blue] (0,0) circle (0.1);}) and $p=250$ periods (\kern0.5pt\protect\tikz[scale=0.3]{\protect\filldraw[color=red] (0,0) circle (0.1);}) with their respective $99\%$ confidence intervals (\protect\tikz\protect\filldraw[fill=blue!40!,draw=none] (0,0) rectangle (0.2,0.2); / \protect\tikz\protect\filldraw[fill=red!40!,draw=none] (0,0) rectangle (0.2,0.2);).}
    \label{fig:frf}
\end{figure}

Next, we use the two sets of obtained FRF measurements to implement FreePC. Since the system is MIMO, we set $E=n_u=2$, and, as outlined in Remark~\ref{rem:FRF-measurements}, we set $\hat{U}_k^e=\bm{e}_e$ and $\hat{Y}^e_k = H(e^{j\hat{\omega}_k})\hat{U}^e_k$ for all $k\in\mathbb{Z}_{[0,M-1]}$ and $e\in\mathcal{E}$. We take the prediction horizon to be $T=50$ and the initial input-output trajectory to be of length $\bar{T}=3n_x=12\geqslant \ell_{\Sigma}$. It is straightforward to verify that $\{\hat{U}^e_{[0,M-1]}\}_{e\in\mathcal{E}}$ are CPE of order $E(2M-1)=158\geqslant T+\bar{T}+n_x=66$. We use the regularization parameters $\lambda_{G}=0.01$ and $\lambda_{\sigma}=1\cdot 10^{4}$. Moreover, we use the stage cost $\ell(y,u)=(y-y_{\mathrm{ss}})^\top Q(y-y_{\mathrm{ss}}) + (u-u_{\mathrm{ss}})^\top R(u-u_{\mathrm{ss}})$ with $Q=I$ and $R=0.01I$, where $y_{\mathrm{ss}}$ and $u_{\mathrm{ss}}$ are the desired steady-state input and output, respectively, satisfying \[\begin{bmatrix}
    A-I & B\\
    C & D
\end{bmatrix}\begin{bmatrix}x_{\mathrm{ss}}\\u_{\mathrm{ss}}
\end{bmatrix} = \begin{bmatrix}
    0\\y_{\mathrm{ss}}
\end{bmatrix}.\] For illustrative purposes, we compute these based on the true system matrices $(A,B,C,D)$. We simulate the system in closed-loop with FreePC for $k\in\mathbb{Z}_{[0,L_{\mathrm{sim}}]}$ with $L_{\mathrm{sim}}=100$, where we use $y_{\mathrm{ss}}=(15,-10)$ for $k\in\mathbb{Z}_{[0,\frac{1}{2}L_{\mathrm{sim}}-1]}$ and $y_{\mathrm{ss}}=(-5,10)$ for $k\in\mathbb{Z}_{[\frac{1}{2}L_{\mathrm{sim}},L_{\mathrm{sim}}]}$. The sets of admissible inputs and outputs are given, respectively, by $\mathbb{U}=[-10,20]^2$ and $\mathbb{Y}=[-18,18]^2$. The resulting input and output trajectories are depicted in Fig.~\ref{fig:freepc-results}. As expected, we observe that the FreePC scheme based on $p=250$ periods performs significantly better than the one with $p=10$. In fact, it closely resembles the simulation results obtained using a model-based  benchmark, i.e., a similarly-tuned MPC scheme (based on the exact model), that are also shown in Fig.~\ref{fig:freepc-results}. This illustrates that measuring longer not only leads to more accurate FRF measurements but also directly improves the achieved performance by FreePC.
\begin{table}[!b]
    \centering
    \caption{Average cost achieved by FreePC for $10^3$ data sets.}
    \label{tab:MonteCarlo}%
    \begin{tabular}{@{}c|cccccc@{}}
        \hline
        & MPC & $p=10$ & $p=25$ & $p=50$ & $p=100$ & $p=250$\\\hline
        Mean $J$ & $641.17$ & $641.36$ & $641.23$ & $641.20$ & $641.19$ & $641.18$
    \end{tabular}
\end{table}
To further investigate the relation between the number of periods and the achieved performance, we perform a Monte Carlo study, in which we implement FreePC using based on data sets generated as discussed before with $y_{\mathrm{ss}}=(15,-10)$ for all $k\in\mathbb{Z}_{[0,L_{\mathrm{sim}}]}$. We consider data sets with $\mathbb{E}(n_kn_k^\top)=\sqrt{0.1}I$, $M=50$ frequencies and $p=\{10,25,50,100,250\}$ periods (after we have again discarded the first $100$ periods to eliminate most transient phenomena), for which we compute the achieved cost \[J = \sum_{\mathclap{k\in\mathbb{Z}_{[0,L_{\mathrm{sim}-1}]}}}\ell(y_k,u_k),\] throughout the duration $L_{\mathrm{sim}}=100$ of the simulation. The average cost and the variance of the cost throughout $1000$ runs are recorded in Table~\ref{tab:MonteCarlo} along with the cost of the model-based benchmark. Since the FRF measurements become more representative of the true system~\eqref{eq:batch-reactor} as the number of periods increases, the average achieved cost and its variance decrease for higher $p$, which can, indeed, be seen in Table~\ref{tab:MonteCarlo}. However, the cost does not reach the benchmark cost, likely due to bias introduced by the regularization terms in~\eqref{eq:FreePC}. Alternative regularization schemes have been proposed (in the context of DeePC) that do not introduce such bias at the cost of being computationally more expensive to implement~\cite{Verheijen2023}. Exploring such alternative regularization schemes in the context of FreePC is an interesting topic for future research.
\begin{figure}[!tb]
    \setlength\fwidth{0.38\textwidth}
    \input{images/freepc-results.tex}
    \caption{Simulation of FreePC using data with $p=10$ (\kern0.5pt\protect\tikz[scale=0.6]{\protect\draw[color=blue, line width=1.2pt] (0,0.1)--(17pt,0.1);\protect\draw[color=white] (0,0)--(0.55,0)}) and $p=250$ (\kern0.5pt\protect\tikz[scale=0.6]{\protect\draw[color=red, line width=1.2pt] (0,0.1)--(0.55,0.1);\protect\draw[color=white] (0,0)--(17pt,0)}) periods and MPC (\kern0.5pt\protect\tikz[scale=0.6]{\protect\draw[color=black, dashed, line width=1.2pt] (0,0.1)--(17pt,0.1);\protect\draw[color=white] (0,0)--(0.55,0)}) with constraints (\kern0.5pt\protect\tikz[scale=0.6]{\protect\draw[color=black, dotted, line width=1.2pt] (0,0.1)--(17pt,0.1);\protect\draw[color=white] (0,0)--(0.55,0)}).}
    \label{fig:freepc-results}
\end{figure}

\section{Conclusions}
\label{sec:conclusions}
In this paper, we presented a novel version of Willems' fundamental lemma based on frequency-domain data, which can be used to fully characterize the input-output behavior of an unknown (possibly unstable) system. In doing so, we also provided a definition of persistence of excitation for frequency-domain data, and we presented extensions of these results to accommodate multiple data sets. To show the potential impact of these new results, we demonstrate their application to data-driven simulation, frequency response evaluation (at other frequencies than available in the off-line data), and data-driven LQR. We also used the frequency-domain fundamental lemma to present a novel data-driven predictive control scheme, called FreePC, that utilizes frequency-domain data \emph{without} the need to convert it into a parametric (state-space) model. In addition to providing a theoretical counterpart to the time-domain fundamental lemma and DeePC, our results bring several practical advantages of working with frequency-domain data into the realm of data-driven control. This includes the ability to collect data in closed loop using a pre-stabilizing controller, important computational benefits, selection of active frequency bands, visual insight, and improved robustness to noise. Moreover, by bridging the gap between recent advances in data-driven analysis and (predictive) control, many existing tools for frequency-domain control/identification can now be exploited as well. Numerical results are presented in which we illustrate some of these benefits of using frequency-domain data.

\section*{Appendix}
\subsection{Proof of Theorem~\ref{thm:fd-WFL}}
Let $(\hat{U}_{[0,M-1]},\hat{X}_{[0,M-1]},\hat{Y}_{[0,M-1]})$ be an input-state-output spectrum of $\Sigma$ in~\eqref{eq:td-system} satisfying Assumption~\ref{asm:ctrb}. Suppose that $\hat{U}_{[0,M-1]}$ is PE of order $L+n_x$. 

\subsubsection{Proof of statement~\ref{item:fd-wfl-1}}
To show that statement~\ref{item:fd-wfl-1} holds, it suffices to show that, for any $b\in\mathbb{C}^{n_x+n_uL},$ \begin{equation}\label{eq:full-row-rank}
    b^\hop\begin{bmatrix}
    F_1(\hat{X}_{[0,M-1]}) & F_1^*(\hat{X}_{[1,M-1]})\\
    F_L(\hat{U}_{[0,M-1]}) & F_L^*(\hat{U}_{[1,M-1]})
\end{bmatrix}=0 \implies b=0.\end{equation} 
To this end, we write $b=(\xi,\eta)$ with $\xi\in\mathbb{C}^{n_x}$ and $\eta\in\mathbb{C}^{n_uL}$ such that the left equality in~\eqref{eq:full-row-rank} holds. Define $b_0\coloneqq(\xi,\eta,0_{n_xn_u})$. From~\eqref{eq:full-row-rank}, we get
\begin{equation}\label{eq:b0Theta}
    b_0^\hop\underbrace{\begin{bmatrix}
        F_1(\hat{X}_{[0,M-1]}) & F_1^*(\hat{X}_{[1,M-1]})\\
        F_{L+n_x}(\hat{U}_{[0,M-1]}) & F_{L+n_x}^*(\hat{U}_{[1,M-1]})
    \end{bmatrix}}_{\eqqcolon \Theta\in\mathbb{C}^{(n_x + (L+n_x)n_u)\times (2M-1)}} = 0.
\end{equation}
Let $\Lambda \coloneqq \diag\{1,e^{j\hat{\omega}_1},\hdots,e^{j\hat{\omega}_{M-1}},e^{-j\hat{\omega}_1},\hdots,e^{-j\hat{\omega}_{M-1}}\}$. It immediately follows from~\eqref{eq:b0Theta} that $b^\hop_0 \Theta\Lambda = 0$. Let $\hat{X}_k^+\coloneqq e^{j\hat{\omega}_k}\hat{X}_k$ and $\hat{U}_k^+\coloneqq e^{j\hat{\omega}_k}\hat{U}_k$ for all $k\in\mathbb{Z}_{[0,M-1]}$. Then,
\begin{multline*}
    0 = b_0^\hop \Theta\Lambda = b_0^\hop\begin{bmatrix}
        F_1(\hat{X}^+_{[0,M-1]}) & F_1^*(\hat{X}^+_{[1,M-1]})\\
        F_{L+n_x}(\hat{U}^+_{[0,M-1]}) & F_{L+n_x}^*(\hat{U}^+_{[1,M-1]})
    \end{bmatrix},\\
    \stackrel{\eqref{eq:fd-system-state}}{=}\begin{bmatrix}
        A^\top\xi\\
        B^\top\xi\\
        \eta\\
        0_{(n_x-1)n_u}
    \end{bmatrix}^\hop \begin{bmatrix}
        F_1(\hat{X}_{[0,M-1]}) & F_1^*(\hat{X}_{[1,M-1]})\\
        F_{L+n_x}(\hat{U}_{[0,M-1]}) & F_{L+n_x}^*(\hat{U}_{[1,M-1]})
    \end{bmatrix},
\end{multline*}
where we exploit the specific structure of $F_1(\hat{X}_{[0,M-1]})$ and $F_{L+n_x}(\hat{U}_{[0,M-1]})$, defined in~\eqref{eq:F_L}.
Repeating this process yields that, for any $m\in\mathbb{Z}_{[0,n_x]}$, 
\begin{equation}\label{eq:bm}
    0 = b_0^\hop \Theta\Lambda^{m} \stackrel{\eqref{eq:fd-system-state}}{=} 
    \underbrace{\begin{bmatrix}\begin{smallmatrix}
        \xi^\hop A^m & \xi^\hop A^{m-1}B & \hdots & \xi^\hop B & \eta^\hop & 0^\top_{(n_x-m)n_u}
    \end{smallmatrix}\end{bmatrix}}_{\eqqcolon b_m^\hop} \Theta.
\end{equation}

Since $\hat{U}_{[0,M-1]}$ is PE of order $L+n_x$, the bottom $(L+n_x)n_u$ rows of $\Theta$ are linearly independent, i.e., $\begin{bmatrix} F_{L+n_x}(\hat{U}_{[0,M-1]}) & F^*_{L+n_x}(\hat{U}_{[1,M-1]})\end{bmatrix}$ has full row rank being $(L+n_x)n_u$. Hence, $\dim\ker\Theta^\hop = n_x+(L+n_x)n_u - \operatorname{rank}\Theta \leqslant n_x$, which means that the $n_x+1$ vectors $b_m$, $m\in\mathbb{Z}_{[0,n_x]}$, as defined in~\eqref{eq:bm}, must be linearly dependent. We can deduce that $\eta=0$ as follows: We write $\eta = (\eta_1,\eta_2,\hdots,\eta_L)$ with $\eta_i\in\mathbb{C}^{ n_u}$ for $i\in\mathbb{Z}_{[1,L]}$. If $\eta_L\neq 0$, then the matrix\begin{equation}\label{eq:upp-triag}\begin{bmatrix} b_0 & \hdots & b_{n_x}\end{bmatrix} = \begin{bmatrix}\begin{smallmatrix}
    \xi & A^\top\xi & \hdots & (A^{n_x})^\top\xi\\
    \eta_1 & B^\top\xi & \hdots & (A^{n_x-1}B)^\top\xi\\
    \eta_2 & \eta_1 & \hdots & (A^{n_x-2}B)^\top\xi\\
    \vdots & \vdots & \ddots & \vdots\\
    \eta_L & \eta_{L-1} & \hdots & \beta_L(\xi,\eta)\\
    0 & \eta_L & \hdots & \beta_{L+1}(\xi,\eta)\\
    \vdots & \vdots & \ddots & \vdots\\
    0 & 0 & \hdots & \eta_{L-1}\\
    0 & 0 & \hdots & \eta_L
\end{smallmatrix}\end{bmatrix},\end{equation}
with \[
\beta_\lambda(\xi,\eta) = \begin{cases}
    (A^{n_x-\lambda}B)^\top\xi,\quad &\text{if }\lambda\leqslant n_x,\\
    \eta_{\lambda-L}, &\text{if }\lambda>n_x,
\end{cases}
\]
would have full column rank, which contradicts the linear dependence of the vectors $b_m$, $m\in\mathbb{Z}_{[0,n_x]}$. Hence, it must hold that $\eta_L=0$. After substituting $\eta_L=0$ in~\eqref{eq:upp-triag}, we see that, if $\eta_{L-1}\neq 0$, the matrix in~\eqref{eq:upp-triag} again would have full column rank, which contradicts with the linear dependence of the vectors $b_m$, $m\in\mathbb{Z}_{[0,n_x]}$. Hence, it must also hold that $\eta_{L-1}=0$. Repeating this process $L$ times yields $\eta=0$. 

By the Cayley-Hamilton theorem, there exist $c_i\in\mathbb{R}$, $i\in\mathbb{Z}_{[0,n_x-1]}$, such that $A^{n_x} = \sum_{n\in\mathbb{Z}_{[0,n_x-1]}} c_nA^n$. It follows that \begin{equation*}
\bar{b} = b_{n_x} - \sum_{\mathclap{n\in\mathbb{Z}_{[0,n_x-1]}}}c_nb_n
= \begin{bmatrix}\begin{smallmatrix}
    0_{n_x\times 1}\\
    ((A^{n_x-1} - c_{n_x-1}A^{n_x-2} - \hdots-c_1I_{n_x})B)^\top\xi\\
    ((A^{n_x-2} - c_{n_x-1}A^{n_x-3} - \hdots -c_2I_{n_x})B)^\top\xi\\
    \vdots\\
    ((A - c_{n_x-1}I_{n_x})B)^\top\xi\\
    B^\top\xi\\
    0_{n_uL\times 1}
\end{smallmatrix}\end{bmatrix}.
\end{equation*}
Since $\bar{b}$ is a linear combination of the vectors $b_m$, $m\in\mathbb{Z}_{[0,n_x]}$, it holds that $\bar{b}^\hop \Theta = 0$. However, due to PE, the bottom $L+n_x$ rows of $\Theta$ are full row rank and, thus, $\bar{b}=0$. It follows that $B^\top\xi = 0$, which we substitute in the row above to find $B^\top (A-c_{n_x}I_{n_x})^\top\xi = B^\top A^\top\xi = 0$. Repeating this process yields $B^\top (A^{n_x-1}-c_{n_x-1}A^{n_x-2}-\hdots - c_1I_{n_x})^\top\xi = B^\top (A^{n_x-1})^\top \xi=0$, which, due to Assumption~\ref{asm:ctrb}, means that $\xi=0$. Thus, we can conclude that $b_0=0$. It follows that~\eqref{eq:full-row-rank} indeed holds, and, thereby, statement~\ref{item:fd-wfl-1} holds.

\subsubsection{Proof of statement~\ref{item:fd-wfl-2}}
By applying~\eqref{eq:td-system} repeatedly, the pair of real-valued sequences $(u_{[0,L-1]},y_{[0,L-1]})$ of length $L$ is an input-output trajectory of $\Sigma$ if and only if there exists $x_0\in\mathbb{R}^{n_x}$ such that \begin{equation}\label{eq:uyx0u}\begin{bmatrix}
    u_{[0,L-1]}\\
    y_{[0,L-1]}
\end{bmatrix} = \begin{bmatrix}
    0 & I\\
    \mathcal{O}_{L} & \mathcal{T}_{L}
\end{bmatrix}\begin{bmatrix}
    x_0\\
    u_{[0,L-1]}
\end{bmatrix},
\end{equation} where $\mathcal{O}_L$ is the observability matrix~\eqref{eq:obs-mat} and
\begin{equation}\label{eq:TL}
    \mathcal{T}_L\coloneqq\begin{bmatrix}\begin{smallmatrix}
        D & 0 & \hdots & 0 & 0\\
        CB & D & \hdots & 0 & 0\\
        \vdots & \vdots & \ddots & \vdots & \vdots\\
        CA^{L-2}B & CA^{L-3}B & \hdots & CB & D
    \end{smallmatrix}\end{bmatrix}.
\end{equation}
Let $\Omega_1 \coloneqq \begin{bmatrix}
    F_1(\hat{X}_{[0,M-1]}) & F^*_1(\hat{X}_{[1,M-1]})\\
    F_L(\hat{U}_{[0,M-1]}) & F^*_L(\hat{U}_{[1,M-1]})
\end{bmatrix}$. By Theorem~\ref{thm:fd-WFL}.\ref{item:fd-wfl-1}, $\Omega_1$ has full row rank. Hence, for any $x_0\in\mathbb{R}^{n_x}$ and $u_{[0,L-1]}\in\mathbb{R}^{n_uL}$, we can express $(x_0,u_{[0,L-1]})$ as $(x_0,u_{[0,L-1]})=\Omega_1 \hat{G}$ for some $\hat{G}\in\mathbb{C}^{2M-1}$. Substitution in~\eqref{eq:uyx0u} shows that $(u_{[0,L-1]},y_{[0,L-1]})$ is an input-output trajectory of $\Sigma$ if and only if there exists $\hat{G}\in\mathbb{C}^{2M-1}$ such that \begin{equation}\label{eq:Omega1G}\begin{bmatrix}
    u_{[0,L-1]}\\
    y_{[0,L-1]}
\end{bmatrix} = \begin{bmatrix}
    0 & I\\
    \mathcal{O}_{L} & \mathcal{T}_{L}
\end{bmatrix}\Omega_1 \hat{G}.
\end{equation}
By~\eqref{eq:fd-system}, $\hat{U}_{[0,M-1]}$, $\hat{X}_{[0,M-1]}$ and $\hat{Y}_{[0,M-1]}$ satisfy \[W_L(e^{j\hat{\omega}_{m}})\otimes \hat{Y}_m = \mathcal{O}_{L}\hat{X}_m + \mathcal{T}_{L}(W_L(e^{j\hat{\omega}_m})\otimes \hat{U}_m),\] for all $m\in\mathbb{Z}_{[0,M-1]}$. Thereby, it follows that
\begin{equation}\label{eq:OLTLF}\begin{bmatrix} F_L(\hat{Y}_{[0,M-1]}) & F_L(\hat{Y}_{[1,M-1]})\end{bmatrix} =\\ \begin{bmatrix}
    \mathcal{O}_L & \mathcal{T}_L
\end{bmatrix}\Omega_1.\end{equation}
Substitution in~\eqref{eq:Omega1G} shows that $(u_{[0,L-1]},y_{[0,L-1]})$ is an input-output trajectory of $\Sigma$ if and only if \begin{equation}\label{eq:FG}\begin{bmatrix}
    u_{[0,L-1]}\\
    y_{[0,L-1]}
\end{bmatrix} = \begin{bmatrix}
    F_L(\hat{U}_{[0,M-1]}) & F^*_L(\hat{U}_{[1,M-1]})\\
    F_L(\hat{Y}_{[0,M-1]}) & F^*_L(\hat{Y}_{[1,M-1]})
\end{bmatrix}\hat{G},\end{equation}
for some $\hat{G}\in\mathbb{C}^{2M-1}$. If $G_0\in\mathbb{R}$ and $G_1\in\mathbb{C}^{2(M-2)}$ satisfy the conditions in Theorem~\ref{thm:fd-WFL}.\ref{item:fd-wfl-2}, then $\hat{G}=(G_0,G_1,G_1^*)$ satisfies~\eqref{eq:FG}. This concludes the proof of \underline{sufficiency} of~\ref{item:fd-wfl-2}.

We proceed by showing the necessity of statement~\ref{item:fd-wfl-2}. To this end, suppose $(u_{[0,L-1]},y_{[0,L-1]})$ is an input-output trajectory of $\Sigma$. Then, as shown above, there exists $\hat{G}\in\mathbb{C}^{2M-1}$ satisfying~\eqref{eq:FG}. Since $(x_0,u_{[0,L-1]})$ is real, it follows, using the particular complex-conjugate structure in~\eqref{eq:FG}, that \[\begin{bmatrix}
    x_0\\
    u_{[0,L-1]}
\end{bmatrix} = \frac{1}{2}\begin{bmatrix}
    F_L(\hat{U}_{[0,M-1]}) & F^*_L(\hat{U}_{[1,M-1]})\\
    F_L(\hat{Y}_{[0,M-1]}) & F^*_L(\hat{Y}_{[1,M-1]})
\end{bmatrix}\begin{bmatrix}
    \hat{G}_0+\hat{G}_0^*\\
    \hat{G}_1+\hat{G}_2^*\\
    \hat{G}_1^*+\hat{G}_2
\end{bmatrix}.\]
Note that $G_0=\frac{1}{2}(\hat{G}_0+\hat{G}_0^*)$ is real and that $G_1=\frac{1}{2}(\hat{G}_1+\hat{G}^*_2) = \frac{1}{2}(\hat{G}_1^*+\hat{G}_2)^*$. Thus, $G_0$ and $G_1$ satisfy the conditions in Theorem~\ref{thm:fd-WFL}.\ref{item:fd-wfl-2}, which completes the proof of \underline{necessity}.

\subsection{Sketch of proof for Theorem~\ref{thm:fd-WFL-multi}}
Let $\{(\hat{U}^e_{[0,M-1]},\hat{X}^e_{[0,M-1]},\hat{Y}^e_{[0,M-1]})\}_{e\in\mathcal{E}}$ be a collection of input-state-output spectra of $\Sigma$ in~\eqref{eq:td-system} satisfying Assumption~\ref{asm:ctrb}. Suppose that $\{\hat{U}^e_{[0,M-1]}\}_{e\in\mathcal{E}}$ is CPE of order $L+n_x$.

\subsubsection{Sketch of proof for statement~\ref{item:fd-wfl-multi-1}}
Similar to the proof of Theorem~\ref{thm:fd-WFL}.\ref{item:fd-wfl-1}, it suffices to show that, for any $b\in\mathbb{C}^{n_x+n_uL}$, \[b^\hop\begin{bmatrix}
    \mathcal{F}_1(\{\hat{X}^e_{[0,M-1]}\}_{e\in\mathcal{E}}) & \mathcal{F}_1^*(\{\hat{X}^e_{[1,M-1]}\}_{e\in\mathcal{E}})\\
    \mathcal{F}_L(\{\hat{U}^e_{[0,M-1]}\}_{e\in\mathcal{E}}) & \mathcal{F}_L^*(\{\hat{U}^e_{[1,M-1]}\}_{e\in\mathcal{E}})
\end{bmatrix}=0\] implies that $b=0$. To show this, we replace \[\Theta\coloneqq \begin{bmatrix}
    \mathcal{F}_1(\{\hat{X}^e_{[0,M-1]}\}_{e\in\mathcal{E}}) & \mathcal{F}_1^*(\{\hat{X}^e_{[1,M-1]}\}_{e\in\mathcal{E}})\\
    \mathcal{F}_{L+n_x}(\{\hat{U}^e_{[0,M-1]}\}_{e\in\mathcal{E}}) & \mathcal{F}_{L+n_x}^*(\{\hat{U}^e_{[1,M-1]}\}_{e\in\mathcal{E}})
\end{bmatrix}\] and $\Lambda\coloneqq\diag\{I_E,\Lambda_1,\Lambda_1^*\}$, with $\Lambda_1=\diag\{e^{j\hat{\omega}_1}I_E,\hdots,e^{j\hat{\omega}_{M-1}}I_E\}$, in the proof of Theorem~\ref{thm:fd-WFL}.\ref{item:fd-wfl-1}, after which the proof holds \textit{mutatis mutandis.}

\subsubsection{Sketch of proof for statement~\ref{item:fd-wfl-multi-2}}
We replace \[\Omega_1\coloneqq\begin{bmatrix}
\mathcal{F}_1(\{\hat{X}^e_{[0,M-1]}\}_{e\in\mathcal{E}}) & \mathcal{F}_1^*(\{\hat{X}^e_{[1,M-1]}\}_{e\in\mathcal{E}})\\
F_L(\{\hat{U}^e_{[0,M-1]}\}_{e\in\mathcal{E}}) & F_L^*(\{\hat{U}^e_{[1,M-1]}\}_{e\in\mathcal{E}})\end{bmatrix},\] in the proof of Theorem~\ref{thm:fd-WFL}.\ref{item:fd-wfl-2}, after which the proof holds \textit{mutatis mutandis}. Some readers may wonder how the different initial conditions for each experiment affect the results. To address this, we point out that, since Definition~\ref{def:input-output-spectra} implicitly assumes that no transient phenomena are present, the initial state only shows up as a phase shift. Due to linearity and time invariance of $\Sigma$, the resulting phase shifts $\phi^e_{[0,M-1]}$, $e\in\mathcal{E}$, affect input, state and output in exactly the same way. Thereby, it holds that $\mathcal{F}_L(\{\hat{Y}^e_{[0,M-1]}\}_{e\in\mathcal{E}})= \mathcal{O}_L\mathcal{F}_1(\{\hat{X}^e_{[0,M-1]}\}_{e\in\mathcal{E}}) + \mathcal{T}_L\mathcal{F}_L(\{\hat{U}^e_{[0,M-1]}\}_{e\in\mathcal{E}})$, with $\mathcal{O}_L$ as in~\eqref{eq:obs-mat} and $\mathcal{T}_L$ as in~\eqref{eq:TL}.

\subsection{Proof of Proposition~\ref{prop:simulation}}
Let $\{(\hat{U}^e_{[0,M-1]},\hat{Y}^e_{[0,M-1]})\}_{e\in\mathcal{E}}$ be a collection of input-output spectra of $\Sigma$ in~\eqref{eq:td-system} satisfying Assumption~\ref{asm:ctrb}. Let $\{\hat{X}^e_{[0,M-1]}\}_{e\in\mathcal{E}}$ be the corresponding state spectra. Suppose that $\{\hat{U}^e_{[0,M-1]}\}_{e\in\mathcal{E}}$ is CPE of order $L+L_0+n_x$. For the past input-output trajectory $(u_{[-L_0,-1]},y_{[-L_0,-1]})$, there exists $x_{-L_0}\in\mathbb{R}^{n_x}$ such that 
\begin{equation}\label{eq:y-past}
    y_{[-L_0,-1]} = \mathcal{O}_{L_0}x_{-L_0} + \mathcal{T}_{L_0}u_{[-L_0,-1]},
\end{equation}
with $\mathcal{O}_L$ as in~\eqref{eq:obs-mat} and $\mathcal{T}_L$ as in~\eqref{eq:TL}. Due to CPE, by Theorem~\ref{thm:fd-WFL-multi}.\ref{item:fd-wfl-multi-1}, the matrix \[\Omega_1\coloneqq\begin{bmatrix}
    \mathcal{F}_{1}(\{\hat{X}^e_{[0,M-1]}\}_{e\in\mathcal{E}}) & \mathcal{F}_{1}^*(\{\hat{X}^e_{[1,M-1]}\}_{e\in\mathcal{E}})\\
    \mathcal{F}_{L_0+L}(\{\hat{U}^e_{[0,M-1]}\}_{e\in\mathcal{E}}) & \mathcal{F}_{L_0+L}^*(\{\hat{U}^e_{[1,M-1]}\}_{e\in\mathcal{E}})
\end{bmatrix}\] has full row rank. It follows that, for any $x_{-L_0}\in\mathbb{R}^{n_x}$ and $u_{[-L_0,L-1]}\in\mathbb{R}^{n_u(L_0+L)}$, there exists a $\hat{G}\in\mathbb{C}^{2M-1}$ such that $(x_{-L_0},u_{[-L_0,L-1]})=\Omega_1\hat{G}$. Hence, using~\eqref{eq:y-past}, there exists $\hat{G}\in\mathbb{C}^{2M-1}$ such that \begin{align}\def\arraystretch{\arraystretchval}&\left[\begin{array}{@{}c@{}}
    u_{[-L_0,-1]}\\
    u_{[0,L-1]}\\\hdashline[2pt/2pt]
    y_{[-L_0,-1]}
\end{array}\right] = \left[\begin{array}{@{}c;{2pt/2pt}cc@{}}
    0 & I & 0\\
    0 & 0 & I\\\hdashline[2pt/2pt]
    \mathcal{O}_{L_0} & \mathcal{T}_{L_0} & 0
\end{array}\right]\Omega_1\hat{G},\label{eq:FG2}\\
&\quad\stackrel{\eqref{eq:fd-system}}{=}\def\arraystretch{\arraystretchval}\left[\begin{array}{@{}cc@{}}
    \mathcal{F}_{L_0+L}(\{\hat{U}^e_{[0,M-1]}\}_{e\in\mathcal{E}}) & \mathcal{F}_{L_0+L}^*(\{\hat{U}^e_{[1,M-1]}\}_{e\in\mathcal{E}})\\\hdashline[2pt/2pt]
    \mathcal{F}_{L_0}(\{\hat{Y}^e_{[0,M-1]}\}_{e\in\mathcal{E}}) & \mathcal{F}_{L_0}^*(\{\hat{Y}^e_{[1,M-1]}\}_{e\in\mathcal{E}})
\end{array}\right]\hat{G}.\nonumber\end{align} 
Since $(u_{[-L_0,-1]},u_{[0,L-1]},y_{[-L_0,-1]})$ is real and using the complex-conjugate structure in~\eqref{eq:FG2}, $G=(G_0,G_1,G_1^*)$ with $G_0=\frac{1}{2}(\hat{G}_0+\hat{G}_0^*)$ and $G_1=\frac{1}{2}(\hat{G}_1+\hat{G}_2^*)$ satisfies \begin{multline*}\def\arraystretch{\arraystretchval}(u_{[-L_0,-1]},u_{[0,L-1]},y_{[-L_0,-1]})=\\ \left[\begin{array}{@{}cc@{}}
    \mathcal{F}_{L_0+L}(\{\hat{U}^e_{[0,M-1]}\}_{e\in\mathcal{E}}) & \mathcal{F}_{L_0+L}^*(\{\hat{U}^e_{[1,M-1]}\}_{e\in\mathcal{E}})\\
    \mathcal{F}_{L_0}(\{\hat{Y}^e_{[0,M-1]}\}_{e\in\mathcal{E}}) & \mathcal{F}_{L_0}^*(\{\hat{Y}^e_{[1,M-1]}\}_{e\in\mathcal{E}})
\end{array}\right]G.\end{multline*} Note that $G_0$ is real. This completes the first part of the proof.

The fact that $(u_{[-L_0,L-1]},y_{[-L_0,L-1]})$ with $y_{[-L_0,L-1]}=\\
    \def\arraystretch{\arraystretchval}\left[\begin{array}{@{}cc@{}}
    \mathcal{F}_{L_0+L}(\{\hat{Y}^e_{[0,M-1]}\}_{e\in\mathcal{E}}) & \mathcal{F}_{L_0+L}^*(\{\hat{Y}^e_{[1,M-1]}\}_{e\in\mathcal{E}})
\end{array}\right]G$ is an input-output trajectory of $\Sigma$ readily follows from Theorem~\ref{thm:fd-WFL-multi}.\ref{item:fd-wfl-multi-2}. If, in addition, $L_0\geqslant \ell_{\Sigma}$, we have that 
\begin{equation}\label{eq:kers}
    \ker \mathcal{O}_{L_0} = \ker \mathcal{O}_{L_0 + m},\text{ for all }m\in\mathbb{Z}_{\geqslant 1}.
\end{equation} 
Note that $x_{-L_0}$ satisfying~\eqref{eq:y-past} may not be unique. In fact, the set of all $x_{-L_0}$ satisfying~\eqref{eq:y-past} is given by $\mathbb{X}_{-L_0}\coloneqq \xi + \ker\mathcal{O}_{L_0}$ with $\xi\in\mathbb{R}^{n_x}$ such that $y_{[-L_0,-1]}=\mathcal{O}_{L_0}\xi + \mathcal{T}_{L_0}u_{[-L_0,-1]}$. Let $\xi\in\mathbb{X}_{-L_0}$. Then, the combined past and future output $y_{[-L_0,L-1]}$ satisfy, for all $\eta\in\ker\mathcal{O}_{L_0}$, $y_{[-L_0,L-1]} = \mathcal{O}_{L_0+L}(\xi+\eta) + \mathcal{T}_{L_0+L}u_{[-L_0,L-1]} \stackrel{\eqref{eq:kers}}{=} \mathcal{O}_{L_0+L}\xi + \mathcal{T}_{L_0+L}u_{[-L_0,L-1]}$. Thereby, $y_{[0,L-1]}$ is unique.

\subsection{Proof of Proposition~\ref{prop:fd-eval}}
Let $\{(\hat{U}^e_{[0,M-1]},\hat{Y}^e_{[0,M-1]})\}_{e\in\mathcal{E}}$ be a collection of input-output spectra of $\Sigma$ satisfying Assumption~\ref{asm:ctrb}. Suppose that $\{\hat{U}^e_{[0,M-1]}\}_{e\in\mathcal{E}}$ is CPE of order $L_0+1+n_x$ with $L_0\in\mathbb{Z}_{\geqslant \ell_{\Sigma}}$. Let $z\in\mathbb{C}$ be a complex frequency that is not an eigenvalue of $A$ and let $U_z\in\mathbb{C}^{n_u}$ be the corresponding input spectrum.

Let $L=L_0+1$. For the (complex-valued) input sequence $u_{[0,L-1]}=W_{L}(z)\otimes U_z$ and initial condition $x_{0} = (zI-A)^{-1}BU_z$, we obtain, for all $k\in\mathbb{Z}_{[0,L-1]}$,
\begin{align}
    &y_k = CA^kx_{0} + C\sum_{j=0}^{k-1}A^{k-1-j}Bu_j + Du_k=\label{eq:y-sol}\\
    &\quad\Big(CA^k(zI-A)^{-1}B + C\sum_{j=0}^{k-1}A^{k-1-j}Bz^j+Dz^k\Big)U_z.\nonumber
\end{align}
If $z=0$, then, since $z$ is not an eigenvalue of $A$, $A$ must be invertible, and we obtain, for all $k\in\mathbb{Z}_{[0,L-1]}$,
\begin{align*}
    y_k &= \begin{cases}
        \left(CA^{k}(-A)^{-1}B + D\right)U_z,&\text{ if }k=0,\\
        \left(CA^{k}(-A)^{-1}B + CA^{k-1}B\right)U_z,&\text{ if }k\in\mathbb{Z}_{[1,L-1]},
    \end{cases}\\
    &= \begin{cases}
        \left(-CA^{k-1}B + D\right)U_z,&\text{ if }k=0,\\
        0,&\text{ if }k\in\mathbb{Z}_{[1,L-1]},
    \end{cases}\\
    &= \left(C(zI-A)^{-1}B + D\right)z^kU_z = H(z)z^kU_z.
\end{align*}
If $z\neq 0$, we obtain, for all $k\in\mathbb{Z}_{[0,L-1]}$,
\begin{align*}
    &y_k = (CA^k(zI-A)^{-1}B + z^{k-1}C\sum_{j=0}^{k-1}(Az^{-1})^{j}B + Dz^k)U_z,\\
    &= \left(CA^k(zI-A)^{-1}B + \right.\\
    &\quad \left.z^{k-1}C(I-(Az^{-1})^k)(I-Az^{-1})^{-1}B + Dz^k\right)U_z,\\
    &= \left(C(zI-A)^{-1}B + D\right)z^kU_z = H(z)z^kU_z,
\end{align*}
where we used the fact that \[\sum_{j=0}^{k-1} (Az^{-1})^{j} = (I-(Az^{-1})^k)(I-Az^{-1})^{-1}.\] Let $Y_z=H(z)U_z$. It follows that both the real-valued pair of time-domain sequences $(\Re(W_{L}(z)\otimes U_z),\Re(W_{L}(z)\otimes Y_z))$ and $(\Im(W_{L}(z)\otimes U_z),\Im(W_{L}(z)\otimes Y_z))$ are input-output trajectories of $\Sigma$ (with initial conditions $x_0=\Re((zI-A)^{-1}BU_z)$ and $x_0=\Im((zI-A)^{-1}BU_z)$, respectively), which, using Theorem~\ref{thm:fd-WFL-multi}.\ref{item:fd-wfl-multi-2} (separately for $(\Re(W_{L}(z)\otimes U_z),\Re(W_{L}(z)\otimes Y_z))$ and $(\Im(W_{L}(z)\otimes U_z),\Im(W_{L}(z)\otimes Y_z))$) holds, if and only if 
\begin{align}
    &\left[\begin{array}{@{}c@{}}
        W_{L}(z)\otimes U_z\\
        W_{L}(z)\otimes Y_z
    \end{array}\right] = \def\arraystretch{\arraystretchval}\left[\begin{array}{@{}c@{}}
        \Re{(W_{L}(z)\otimes U_z)}+ j\Im{(W_{L}(z)\otimes U_z)}\\
        \Re{(W_{L}(z)\otimes Y_z)}+j\Im{(W_{L}(z)\otimes Y_z)}
    \end{array}\right] = \nonumber\\
    \def\arraystretch{\arraystretchval}&\quad\left[\begin{array}{@{}c;{2pt/2pt}c@{}}
        \mathcal{F}_{L}(\{\hat{U}^e_{[0,M-1]}\}_{e\in\mathcal{E}}) & \mathcal{F}^*_{L}(\{\hat{U}^e_{[1,M-1]}\}_{e\in\mathcal{E}})\\\hdashline[2pt/2pt]
        \mathcal{F}_{L}(\{\hat{Y}^e_{[0,M-1]}\}_{e\in\mathcal{E}}) & \mathcal{F}^*_{L}(\{\hat{Y}^e_{[1,M-1]}\}_{e\in\mathcal{E}})
    \end{array}\right]G,\label{eq:fd-eval-rewritten}
\end{align}
for some $G_{0,r},G_{0,i}\in\mathbb{R}^{E}$ and $G_{1,r},G_{1,i}\in\mathbb{R}^{E(M-1)}$, where \begin{equation*}
    G = \left[\begin{array}{@{}c@{}}
        G_{0,r} + jG_{0,i}\\
        G_{1,r} + jG_{1,i}\\\hdashline[2pt/2pt]
        G_{1,r}^* + jG_{1,i}^*
    \end{array}\right]\label{eq:G-sum}=\underbrace{\begin{bmatrix}
        I & jI & 0 & 0 & 0 & 0\\
        0 & 0 & I & jI & jI & -I\\
        0 & 0 & I & -jI & jI & I
    \end{bmatrix}}_{\eqqcolon T_G\in\mathbb{C}^{E(2M-1)\times 2E(2M-1)}}\mathcal{G},
\end{equation*}
where $\mathcal{G} = (G_{0,r},G_{0,i},\Re G_{1,r},\Im G_{1,r},\Re G_{1,i},\Im G_{1,i})$. Since $\operatorname{rank}T_G = E(2M-1)$, it follows that $(W_{L}(z)\otimes U_z,W_{L}(z)\otimes Y_z)$ is a (complex-valued) input-output trajectory of $\Sigma$ (with $x_0 =(zI-A)^{-1}BU_z$), if and only if there exists $G\in\mathbb{C}^{E(2M-1)}$ satisfying~\eqref{eq:fd-eval-rewritten}. Note that~\eqref{eq:fd-eval-rewritten} is~\eqref{eq:fd-eval} rewritten. 


We complete the proof by showing that $Y_z$ is unique. To this end, note that, since $L_0\geqslant \ell_{\Sigma}$, we have that~\eqref{eq:kers} holds. We can, equivalently, write~\eqref{eq:y-sol} as \begin{equation*}\label{eq:OL0TL0}
    W_{L_0}(z)\otimes Y_z  = y_{[0,L_0-1]} = \mathcal{O}_{L_0}x_0 + \mathcal{T}_{L_0}(W_{L_0}(z)\otimes U_z).
\end{equation*}
Note that $x_0$ is not necessarily unique, however, it must belong to the set $\mathbb{X}_{0}\coloneqq \xi + \ker\mathcal{O}_{L_0}$ with $\xi\in\mathbb{C}^{n_x}$ such that $W_{L_0}(z)\otimes Y_z=\mathcal{O}_{L_0}\xi + \mathcal{T}_{L_0}(W_{L_0}(z)\otimes U_z)$. Let $\xi\in\mathbb{X}_{0}$. Then, $y_{[0,L_0]}=W_{L_0+1}(z)\otimes Y_z$ satisfies
\begin{align*}
    W_{L_0+1}\otimes Y_z &= \mathcal{O}_{L_0+1}(\xi+\eta) + \mathcal{T}_{L_0+1}(W_{L_0+1}(z)\otimes U_z),\\
    &\stackrel{\eqref{eq:kers}}{=} \mathcal{O}_{L_0+1}\xi + \mathcal{T}_{L_0+1}(W_{L_0+1}(z)\otimes U_z),
\end{align*}
for all $\eta\in\ker\mathcal{O}_{L_0}$. Thus, $y_{[0,L_0]}$ satisfying~\eqref{eq:y-sol} is unique. It follows that $W_{L_0+1}(z)\otimes Y_z$ in~\eqref{eq:fd-eval-rewritten} and, thereby, $Y_z$ (i.e., the first $n_y$ entries of $W_{L_0+1}(z)\otimes Y_z$) satisfying~\eqref{eq:fd-eval} are unique.

\subsection{Proof of Proposition~\ref{prop:lqr}}
Let $(\hat{U}_{[0,M-1]},\hat{Y}_{[0,M-1]})$ be an input-output spectrum of $\Sigma$ in~\eqref{eq:td-system} satisfying Assumption~\ref{asm:ctrb}. Suppose that $\hat{U}_{[0,M-1]}$ is PE of order $n_x+1$ and that $\rank \begin{bmatrix} Q & \lambda I-A\end{bmatrix}= n_x$ for all $\lambda\in\mathbb{C}$ with $|\lambda|=1$. Then, by~\cite[Theorem 23]{vanWaarde2020c}, the optimal control law that achieves $\lim_{k\rightarrow\infty}x_k=0$ and minimizes the cost~\eqref{eq:lqr-cost} is~\eqref{eq:lqr} with \begin{equation}\label{eq:opt-fb-gain}
    \bar{K}=-(R+B^\top \bar{P}B)^{-1}B^\top\bar{P}A,
\end{equation}
where $\bar{P}$ is the \emph{largest} real symmetric solution to the DARE \begin{equation}\label{eq:dare}A^\top PA-P - A^\top PB(R+B^\top PB)^{-1}B^\top PA+Q=0,\end{equation} and satisfies $\bar{P}\succcurlyeq 0$. There is only one largest $\bar{P}$ in the sense that any other solution $\tilde{P}$ to~\eqref{eq:dare} with $\tilde{P}\succcurlyeq P$ for all $P$ satisfying~\eqref{eq:dare} must satisfy $\tilde{P}=\bar{P}$~\cite{vanWaarde2020c}. It follows that 
\begin{equation}\label{eq:dare-rewritten}
    (A+B\bar{K})^\top\bar{P}(A+B\bar{K}) - \bar{P} + \bar{K}^\top R\bar{K} + Q = 0.
\end{equation}

Next, we show that $\bar{P}$ is a solution to~\eqref{eq:lqr-opt}. To this end, let $P=P^\top\in\mathbb{R}^{n_x\times n_x}$ be any symmetric matrix that satisfies $\Delta^\hop\Pi(P)\Delta\succcurlyeq 0$. Since $X_1 = AX_0+BU$, it holds that \begin{equation}\label{eq:X0Udare}\Delta^\hop \Pi(P)\Delta = (\star)^\hop\begin{bmatrix}
    Q-P+A^\top PA & A^\top PB\\
    B^\top PA & R+B^\top PB
\end{bmatrix}\begin{bmatrix}
    X_0\\
    U
\end{bmatrix},\end{equation} where, due to Theorem~\ref{thm:fd-WFL-multi}.\ref{item:fd-wfl-multi-1}, the matrix $\begin{bmatrix} X_0^\top & U^\top\end{bmatrix}^\top$ has full row rank such that $\im\begin{bmatrix}
    I\\
    \bar{K}
\end{bmatrix}\subset\im\begin{bmatrix}X_0\\ U\end{bmatrix}=\mathbb{C}^{n_u+n_x}$. Since $\Delta^\hop\Pi(P)\Delta\succcurlyeq0$, it follows that
\begin{multline*}
    0\preccurlyeq (\star)^\top\begin{bmatrix}
        Q -P+A^\top PA & A^\top PB\\
        B^\top PA & R+B^\top PB
    \end{bmatrix}\begin{bmatrix}
        I\\
        \bar{K}
    \end{bmatrix}=\\
    (A+B\bar{K})^\top P(A+B\bar{K}) -P + \bar{K}^\top R\bar{K} + Q.
\end{multline*}
Subtracting~\eqref{eq:dare-rewritten} from this yields $(\bar{P}-P)-(A+B\bar{K})^\top(\bar{P}-P)(A+B\bar{K}) \succcurlyeq 0$, which, since $A+B\bar{K}$ is Schur, implies that $\bar{P}\succcurlyeq P$. We conclude that $\bar{P}\succcurlyeq P$ for any $P$ that satisfies $\Delta^\hop\Pi(P)\Delta\succcurlyeq 0$, whereby $\trace\bar{P}\geqslant \trace P$. In addition, $\bar{P}$ satisfies $\Delta^\hop \Pi(\bar{P})\Delta=0\succcurlyeq 0$. We conclude that $\bar{P}$ is a solution to~\eqref{eq:lqr-opt}.

To show that $\bar{P}$ is the unique solution to~\eqref{eq:lqr-opt}, let $\tilde{P}$ be another solution to~\eqref{eq:lqr-opt}. Then, $\tilde{P}=\tilde{P}^\top\succcurlyeq 0$, $\trace \tilde{P}=\trace \bar{P}$, and $\Delta^\hop \Pi(\tilde{P})\Delta\succcurlyeq 0$. It follows, as we have shown above, that $\bar{P}\succcurlyeq \tilde{P}$. Thereby, all eigenvalues of $\bar{P}-\tilde{P}$, denoted $\{\lambda_i(\bar{P}-\tilde{P})\}_{i=1}^{n_x}$, are non-negative, i.e., $\bar{P}-\tilde{P}\succcurlyeq 0$. In addition, we know that $0 = \trace\bar{P}-\tilde{P} = \sum_{i=1}^{n_x}\lambda_i(\bar{P}-\tilde{P})$, which, combined with $\bar{P}-\tilde{P}\succcurlyeq 0$, implies that all eigenvalues of $\bar{P}-\tilde{P}$ are zero. Since $\bar{P}-\tilde{P}$ is symmetric, it follows that $\bar{P}-\tilde{P}=0$, and, hence, $\bar{P}$ is the unique solution to~\eqref{eq:lqr-opt}.

It remains to show that $\bar{K}$, i.e., the optimal control law~\eqref{eq:lqr} that minimizes~\eqref{eq:lqr-cost} and achieves $\lim_{k\rightarrow\infty}x_k=0$, can be computed as~\eqref{eq:K-opt}. To this end, observe that we can express 
\begin{align*}
    &\Delta^\hop\Pi(\bar{P})\Delta\stackrel{\eqref{eq:X0Udare}}{=}(\star)^\hop \begin{bmatrix}
        A^\top\bar{P}A-\bar{P}+Q & A^\top\bar{P}B\\
        B^\top\bar{P}A & R+B^\top\bar{P}B
    \end{bmatrix}\begin{bmatrix} X_0\\ U\end{bmatrix}\\
    &\stackrel{\eqref{eq:dare}}{=} (\star)^\hop\begin{bmatrix}
        (R+B^\top\bar{P}B)^{-1} & I\\
        I & R+B^\top\bar{P}B
    \end{bmatrix}\begin{bmatrix} B^\top\bar{P}AX_0\\ U\end{bmatrix}.
\end{align*}
Let $\bar{K}$ be as in~\eqref{eq:opt-fb-gain}. It follows that
\begin{align}
    &\Delta^\hop \Pi\Delta = (\star)^\hop\begin{bmatrix}
        R+B^\top\bar{P}B & -(R+B^\top\bar{P}B)\\
        -(R+B^\top\bar{P}B) & R+B^\top\bar{P}B
    \end{bmatrix}\begin{bmatrix} \bar{K}X_0\\ U\end{bmatrix}\nonumber\\
    &=(U-\bar{K}X_0)^\hop(R+B^\top\bar{P}B)(U-\bar{K}X_0),\label{eq:DeltaPsiDelta-rewritten}
\end{align}
Due to Theorem~\ref{thm:fd-WFL-multi}.\ref{item:fd-wfl-multi-1}, there exists $\begin{bmatrix}
    X_0\\
    U
\end{bmatrix}\Gamma = \begin{bmatrix}I\\ \bar{K}\end{bmatrix}$, for some $\Gamma\in\mathbb{C}^{2M-1\times n_x}$. This implies that $\Gamma$ is a right inverse of $X_0$, whereby $\Gamma X_0\Gamma = \Gamma$ and, since $\Delta^\hop\Pi(\bar{P})\Delta=0$, \[\Delta^\hop\Pi(\bar{P})\Delta\Gamma = (U-\bar{K}X_0)^\hop (R+B^\top \bar{P}B)(U\Gamma-U\Gamma X_0\Gamma)=0.\] 
This shows that there exists a right inverse $X_0^\dagger$ of $X_0$ such that $\Delta^\hop \Pi(\bar{P})\Delta X_0^\dagger = 0$. If $X_0^\dagger$ is a right inverse of $X_0$ satisfying $\Delta^\hop \Pi(\bar{P})\Delta X_0^\dagger=0$, then it follows from~\eqref{eq:DeltaPsiDelta-rewritten}, due to $R+B^\top\bar{P}B\succ 0$, that $(U-\bar{K}X_0)X_0^\dagger=0$. We conclude that the optimal feedback gain satisfies $K=\bar{K}=UX_0^\dagger$.

\subsection{Proof of Theorem~\ref{thm:equivalence}}
Let $\{(\hat{U}^e_{[0,M-1]},\hat{Y}^e_{[0,M-1]})\}_{e\in\mathcal{E}}$ and $(\hat{u}_{[0,N-1]},\hat{y}_{[0,N-1]})$ be, respectively, a collection of input-output spectra and a trajectory of $\Sigma$ in~\eqref{eq:td-system} satisfying Assumption~\ref{asm:ctrb}. Suppose that $\{\hat{U}^e_{[0,M-1]}\}_{e\in\mathcal{E}}$ and $\hat{u}_{[0,M-1]}$ are, respectively, CPE and PE of order $\bar{T}+T+n_x$. For $\lambda_g=\lambda_G=0$,~\eqref{eq:FreePC} and~\eqref{eq:DeePC} only differ in the precise matrix-vector product used in the prediction model to describe the input-output trajectory $(u_{[k-\bar{T},k-1]},u_{[0,T-1],k},y_{[k-\bar{T},k-1]}+\sigma_k,y_{[0,T-1],k})$.
By Lemma~\ref{lem:td-WFL}.\ref{item:td-wfl-2} and Theorem~\ref{thm:fd-WFL-multi}.\ref{item:fd-wfl-multi-2}, both prediction models describe the same space of input-output trajectories so that~\eqref{eq:FreePC} and~\eqref{eq:DeePC} admit the same (set of) optimal control action(s).

\section*{References}
\bibliographystyle{IEEEtran}
\bibliography{../../../../../bibliography/references.bib}

\begin{thebibliography}{10}
\providecommand{\url}[1]{#1}
\csname url@rmstyle\endcsname
\providecommand{\newblock}{\relax}
\providecommand{\bibinfo}[2]{#2}
\providecommand\BIBentrySTDinterwordspacing{\spaceskip=0pt\relax}
\providecommand\BIBentryALTinterwordstretchfactor{4}
\providecommand\BIBentryALTinterwordspacing{\spaceskip=\fontdimen2\font plus
\BIBentryALTinterwordstretchfactor\fontdimen3\font minus \fontdimen4\font\relax}
\providecommand\BIBforeignlanguage[2]{{%
\expandafter\ifx\csname l@#1\endcsname\relax
\typeout{** WARNING: IEEEtran.bst: No hyphenation pattern has been}%
\typeout{** loaded for the language `#1'. Using the pattern for}%
\typeout{** the default language instead.}%
\else
\language=\csname l@#1\endcsname
\fi
#2}}

\bibitem{Berberich2020}
J.~{Berberich} and F.~{Allg\"{o}wer}, ``A trajectory-based framework for data-driven system analysis and control,'' in \emph{Eur. Control Conf.}, 2020, pp. 1365--1370.

\bibitem{Verhoek2021}
C.~{Verhoek}, R.~{T\'{o}th}, S.~{Haesaert}, and A.~{Koch}, ``Fundamental lemma for data-driven analysis of linear parameter-varying systems,'' in \emph{60th IEEE Conf. Decis. Control}, 2021, pp. 5040--5046.

\bibitem{Willems2005}
J.~C. {Willems}, P.~{Rapisarda}, I.~{Markovsky}, and B.~L.~M. {De Moor}, ``A note on persistency of excitation,'' \emph{Syst. Control Lett.}, vol.~54, no.~4, pp. 325--329, 2005.

\bibitem{vanWaarde2020}
H.~J. {van Waarde}, C.~{De Persis}, M.~K. {Camlibel}, and P.~{Tesi}, ``Willems' fundamental lemma for state-space systems and its extension to multiple datasets,'' \emph{IEEE Control Syst. Lett.}, vol.~4, no.~3, pp. 602--607, 2020.

\bibitem{Verhaegen1992}
M.~{Verhaegen} and P.~{Dewilde}, ``Subspace model identification {Part} 1. {The} output-error state-space model identification class of algorithms,'' \emph{Int. J. Control}, vol.~56, no.~5, pp. 1187--1210, 1992.

\bibitem{Molodchyk2024}
O.~{Molodchyk} and T.~{Faulwasser}, ``Exploring the links between the fundamental lemma and kernel regression,'' \emph{IEEE Control Syst. Lett.}, vol.~8, pp. 2045--2050, 2024.

\bibitem{Schmitz2022}
P.~{Schmitz}, T.~{Faulwasser}, and K.~{Worthmann}, ``{Willems}' fundamental lemma for linear descriptor systems and its use for data-driven output-feedback {MPC},'' \emph{IEEE Control Syst. Lett.}, vol.~6, pp. 2443--2448, 2022.

\bibitem{Faulwasser2023}
T.~{Faulwasser}, R.~{Ou}, G.~{Pan}, P.~{Schmitz}, and K.~{Worthmann}, ``Behavioral theory for stochastic systems? {A} data-driven journey from {Willems} to {Wiener} and back again,'' \emph{Annu. Rev. Control}, vol.~55, pp. 92--117, 2023.

\bibitem{Pan2022}
G.~{Pan}, R.~{Ou}, and T.~{Faulwasser}, ``On a stochastic fundamental lemma and its use for data-driven optimal control,'' \emph{IEEE Trans. Autom. Control}, 2022.

\bibitem{Rapisarda2023}
R.~{Rapisarda}, M.~K. {Çamlibel}, and H.~J. {van Waarde}, ``A ``fundamental lemma'' for continuous-time systems, with applications to data-driven simulation,'' \emph{Syst. Control Lett.}, vol. 179, p. 105603, 2023.

\bibitem{Coulson2019}
J.~{Coulson}, J.~{Lygeros}, and F.~{D\"{o}rfler}, ``Data-enabled predictive control: {In} the shallows of the {DeePC},'' in \emph{Eur. Control Conf.}, 2019, pp. 307--312.

\bibitem{Verheijen2023}
P.~C.~N. {Verheijen}, V.~{Breschi}, and M.~{Lazar}, ``Handbook of linear data-driven predictive control: {Theory}, implementation and design,'' \emph{Annu. Rev. Control}, vol.~56, p. 100914, 2023.

\bibitem{Verhoek2021b}
C.~{Verhoek}, H.~S. {Abbas}, R.~{T\'{o}th}, and S.~{Haesaert}, ``Data-driven predictive control for linear parameter-varying systems,'' \emph{IFAC-PapersOnLine}, vol.~54, no.~8, pp. 101--108, 2021.

\bibitem{Bilgic2024}
D.~{Bilgic}, A.~{Harding}, and T.~{Faulwasser}, ``Data-driven predictive control of bilinear {HVAC} dynamics---{An} experimental case study,'' \emph{IEEE Control Syst. Lett.}, vol.~8, pp. 3009--3014, 2024.

\bibitem{Berberich2022b}
J.~{Berberich}, J.~{K\"{o}hler}, M.~A. {M\"{u}ller}, and F.~{Allg\"{o}wer}, ``Linear tracking {MPC} for nonlinear systems--{Part II: The} data-driven case,'' \emph{IEEE Trans. Autom. Control}, vol.~67, no.~9, pp. 4406--4421, 2022.

\bibitem{Alsati2023}
M.~{Alsati}, V.~G. {Lopez}, J.~{Berberich}, F.~{Allg\"{o}wer}, and M.~A. {M\"{u}ller}, ``Data-driven nonlinear predictive control for feedback linearizable systems,'' \emph{IFAC-PapersOnLine}, vol.~56, no.~2, pp. 617--624, 2023.

\bibitem{Lazar2023}
M.~{Lazar}, ``Basis functions nonlinear data-enabled predictive control: {Consistent} and computationally efficient formulations,'' 2023, preprint: \url{https://arxiv.org/abs/2311.05360}.

\bibitem{Azarbahram2024}
A.~{Azarbahram}, M.~{Al Khatib}, V.~K. {Mishra}, and N.~{Bajcinca}, ``Data-driven predictive control for a class of nonlinear systems with polynomial terms,'' \emph{IFAC-PapersOnLine}, vol.~58, no.~21, pp. 226--231, 2024.

\bibitem{Pan2023}
G.~{Pan}, R.~{Ou}, and T.~{Faulwasser}, ``Towards data-driven stochastic predictive control,'' \emph{Int. J. Robust Nonlinear Control}, pp. 1--23, 2023.

\bibitem{Breschi2023}
V.~{Breschi}, A.~{Chiuso}, and S.~{Formentin}, ``Data-driven predictive control in a stochastic setting: {A} unified framework,'' \emph{Automatica}, vol. 152, p. 110961, 2023.

\bibitem{Franklin2010}
G.~F. {Franklin}, J.~D. {Powell}, and A.~{Emami-Naeini}, \emph{Feedback control of dynamic systems}, 6th~ed.\hskip 1em plus 0.5em minus 0.4em\relax Pearson, 2010.

\bibitem{Skogestad2005}
S.~{Skogestad} and I.~{Postlethwaite}, \emph{Multivariable feedback control: {Analysis} and design}, 2nd~ed.\hskip 1em plus 0.5em minus 0.4em\relax Wiley, 2005.

\bibitem{Soderstrom1989}
T.~{S\"{o}derstr\"{o}m} and P.~{Stoica}, \emph{System identification}.\hskip 1em plus 0.5em minus 0.4em\relax Prentice Hall, 1989.

\bibitem{Pintelon2012}
R.~{Pintelon} and J.~{Schoukens}, \emph{System Identification: {A} Frequency Domain Approach}.\hskip 1em plus 0.5em minus 0.4em\relax Wiley, 2012.

\bibitem{Schoukens2004}
J.~{Schoukens}, R.~{Pintelon}, and Y.~{Rolain}, ``Time domain identification, frequency domain identification. equivalences! differences?'' in \emph{Amer. Control Conf.}, 2004, pp. 661--666.

\bibitem{Ferizbegovic2021}
M.~{Ferizbegovic}, H.~{Hjalmarsson}, P.~{Mattsson}, and T.~B. {Sch\"{o}n}, ``{Willems}' fundamental lemma based on second-order moments,'' in \emph{60th IEEE Conf. Decis. Control}, 2021, pp. 396--401.

\bibitem{Markovsky2024}
I.~{Markovsky} and H.~{Ossareh}, ``Finite-data nonparametric frequency response evaluation without leakage,'' \emph{Automatica}, vol. 159, p. 111351, 2024.

\bibitem{Ozkan2012}
L.~{\"{O}zkan}, J.~{Meijs}, and A.~C.~P.~M. {Backx}, ``A frequency domain approach for {MPC} tuning,'' \emph{Comput. Aided Chem. Eng.}, vol.~31, pp. 1632--1636, 2012.

\bibitem{Burgos2014}
J.~G. {Burgos}, C.~A. {L\'{o}pez Mart\'{i}nez}, R.~{van de Molengraft}, and M.~{Steinbuch}, ``Frequency domain tuning method for unconstrained linear output feedback model predictive control,'' in \emph{Proc. 19th IFAC World Congress}, 2014, pp. 7455--7460.

\bibitem{Shah2013}
G.~{Shah} and S.~{Engell}, ``Multivariable {MPC} design based on a frequency response approximation approach,'' in \emph{Eur. Control Conf}, 2013, pp. 13--18.

\bibitem{Sathyanarayanan2023}
K.~K. {Sathyanarayanan}, G.~{Pan}, and T.~{Faulwasser}, ``Towards data-driven predictive control using wavelets,'' \emph{IFAC-PapersOnLine}, vol.~56, no.~2, pp. 632--637, 2023.

\bibitem{Meijer2024-nmpc}
T.~J. {Meijer}, S.~A.~N. {Nouwens}, K.~J.~A. {Scheres}, V.~S. {Dolk}, and W.~P.~M.~H. {Heemels}, ``Frequency-domain data-driven predictive control,'' \emph{IFAC-PapersOnLine}, vol.~58, no.~18, pp. 86--91, 2024.

\bibitem{DePersis2020}
C.~{De Persis} and P.~{Tesi}, ``Formulas for data-driven control: {Stabilization}, optimality, and robustness,'' \emph{IEEE Trans. Autom. Control}, vol.~65, no.~3, pp. 909--924, 2020.

\bibitem{DePersis2021}
------, ``Low-complexity learning of linear quadratic regulators from noisy data,'' \emph{Automatica}, vol. 128, p. 109548, 2021.

\bibitem{Markovsky2005}
I.~{Markovsky}, J.~C. {Willems}, P.~{Rapisarda}, and B.~L.~M. {De Moor}, ``Algorithms for deterministic balanced subspace identification,'' \emph{Automatica}, vol.~41, pp. 755--766, 2005.

\bibitem{Markovsky2006}
I.~{Markovsky}, J.~C. {Willems}, S.~{Van Huffel}, and B.~{De Moor}, \emph{Exact and approximate modeling of linear systems: {A} behavioral approach}.\hskip 1em plus 0.5em minus 0.4em\relax SIAM, 2006.

\bibitem{Berberich2021}
J.~{Berberich}, J.~{K\"{o}hler}, M.~A. {M\"{u}ller}, and F.~{Allg\"{o}wer}, ``Data-driven model predictive control with stability and robustness guarantees,'' \emph{IEEE Trans. Autom. Control}, vol.~66, no.~4, 2021.

\bibitem{McKelvey1996}
T.~{McKelvey}, H.~{Akcay}, and L.~{Ljung}, ``Subspace-based multivariable system identification from frequency response data,'' \emph{IEEE Trans. Autom. Control}, vol.~41, no.~7, pp. 960--979, 1996.

\bibitem{Overschee1996}
P.~{Van Overschee} and B.~{De Moor}, ``Continuous-time frequency domain subspace system identification,'' \emph{Signal Process.}, vol.~52, pp. 179--194, 1996.

\bibitem{Cauberghe2006}
B.~{Cauberghe}, P.~{Guillaume}, R.~{Pintelon}, and P.~{Verboven}, ``Frequency-domain subspace identification using {FRF} data from arbitrary signals,'' \emph{J. Sound Vib.}, vol. 290, pp. 555--571, 2006.

\bibitem{Meijer2025-cdc-arxiv}
T.~J. {Meijer}, M.~{Wind}, V.~S. {Dolk}, and W.~P.~M.~H. {Heemels}, ``Leveraging non-steady-state frequency-domain data in {Willems}' fundamental lemma,'' 2025, preprint: \url{https://arxiv.org/abs/2504.06403}.

\bibitem{Meijer2024-ecc-arxiv}
T.~J. {Meijer}, S.~A.~N. {Nouwens}, V.~S. {Dolk}, and W.~P.~M.~H. {Heemels}, ``A frequency-domain version of {Willems'} fundamental lemma,'' 2023, preprint: \url{https://arxiv.org/abs/2311.15284}.

\bibitem{Proakis1996}
J.~G. {Proakis} and D.~G. {Manolakis}, \emph{Digital Signal Processing: {Principles}, Algorithms and Applications}.\hskip 1em plus 0.5em minus 0.4em\relax Prentice Hall, 1996.

\bibitem{Raff2006}
T.~{Raff}, S.~{Huber}, Z.~K. {Nagy}, and F.~{Allgöwer}, ``Nonlinear model predictive control of a four tank system: {An} experimental stability study,'' in \emph{Proc. Int. Conf. Control Appl.}, 2006, pp. 237--242.

\bibitem{Markovsky2008}
I.~{Markovsky} and P.~{Rapisarda}, ``Data-driven simulation and control,'' \emph{Int. J. Control}, vol.~81, no.~12, pp. 1946--1959, 2008.

\bibitem{vanWaarde2020c}
H.~J. {van Waarde}, J.~{Eising}, H.~L. {Trentelmann}, and M.~K. {Camlibel}, ``Data informativity: {A} new perspective on data-driven analysis and control,'' \emph{IEEE Trans. Autom. Control}, vol.~65, no.~11, pp. 4753--4768, 2020.

\bibitem{Walsh2001}
G.~C. {Walsh} and H.~{Ye}, ``Scheduling of networked control systems,'' \emph{IEEE Control Syst. Mag.}, vol.~21, no.~1, pp. 57--65, 2001.

\bibitem{vandenHof1996}
P.~M.~J. {van den Hof} and R.~A. {de Callafon}, ``Multivariable closed-loop identification: {From} indirect identification to dual-{Youla} parametrization,'' in \emph{35th IEEE Conf. Decis. Control}, 1996, pp. 1397--1402.

\bibitem{vanBerkel2015}
M.~{van Berkel}, ``Estimation of heat transport coefficients in fusion plasmas,'' Ph.D. dissertation, Technische Universiteit Eindhoven, 2015.

\end{thebibliography}

\vspace*{-1cm}
\begin{IEEEbiography}[{\includegraphics[width=1in,height=1.25in,clip,keepaspectratio]{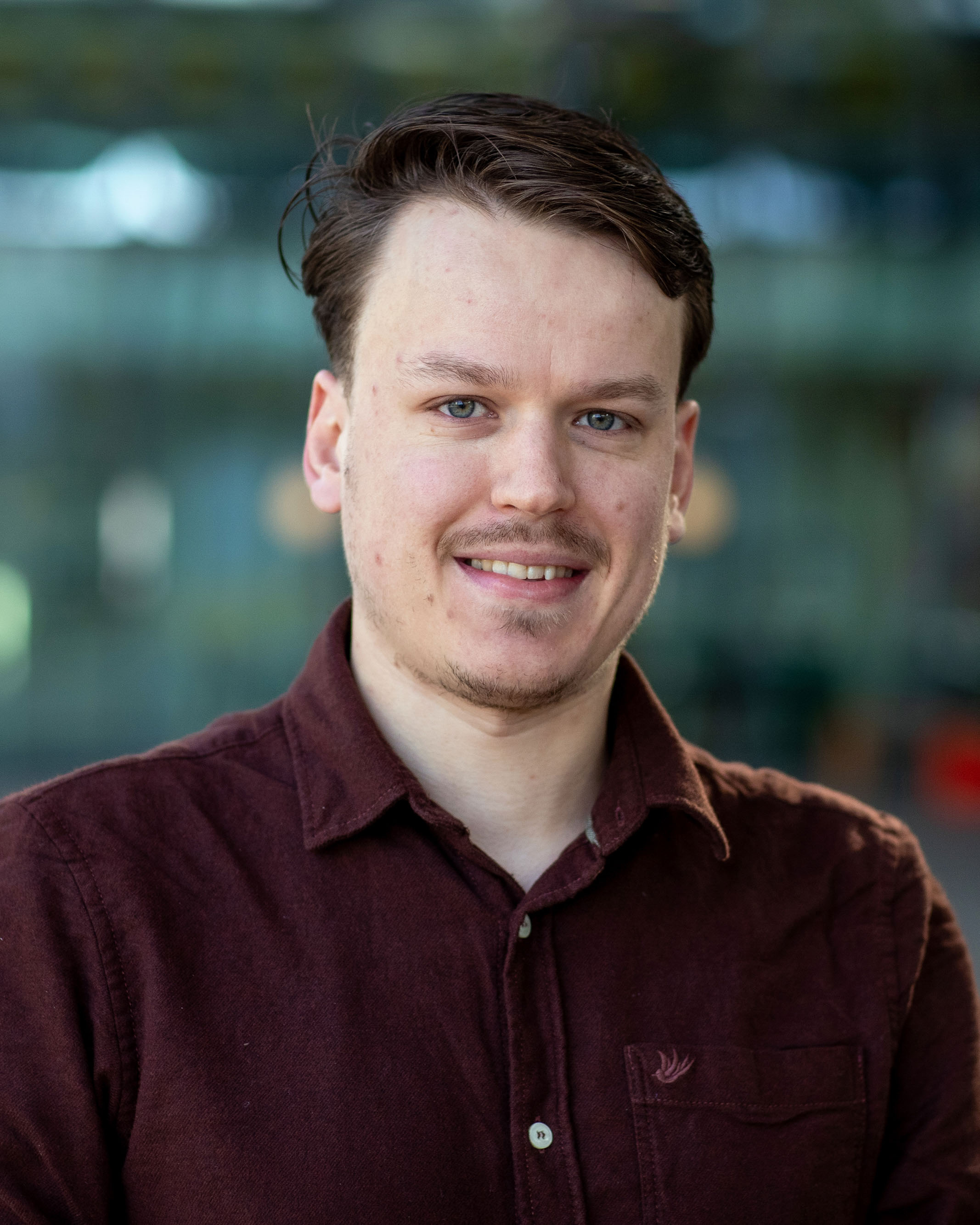}}]{Tomas J. Meijer} (Member, IEEE) received the M.Sc. degree (cum laude) in electrical engineering, in 2018, and the Ph.D. degree in mechanical engineering, in 2024, from the Eindhoven University of Technology (TU/e), the Netherlands, where he is currently a Postdoctoral Researcher. His research interests include data-driven control, model predictive control, parameter-varying systems and distributed-parameter systems.
\end{IEEEbiography}\vspace*{-.5cm}
\begin{IEEEbiography}[{\includegraphics[width=1in,height=1.25in,clip,keepaspectratio]{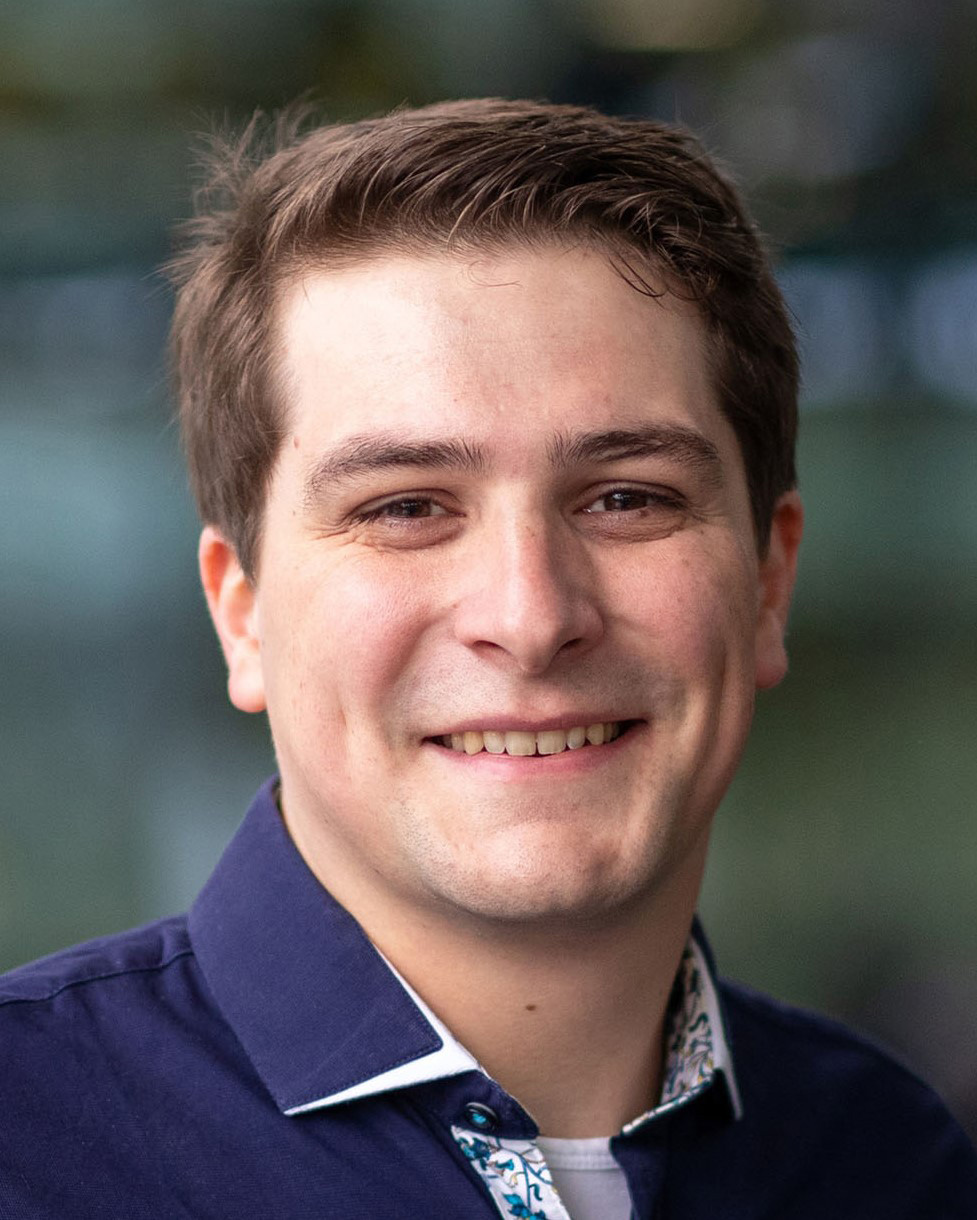}}]{Koen J.A. Scheres} (Member, IEEE) received the M.Sc. degree in mechanical engineering and the Ph.D. degree in control theory from the Eindhoven University of Technology (TU/e), the Netherlands, in 2020 and 2024, respectively. His research interests include data-driven control, hybrid dynamical systems, event-triggered control and neuromorphic control.
\end{IEEEbiography}\vspace*{-.5cm}
\begin{IEEEbiography}[{\includegraphics[width=1in,height=1.25in,clip,keepaspectratio]{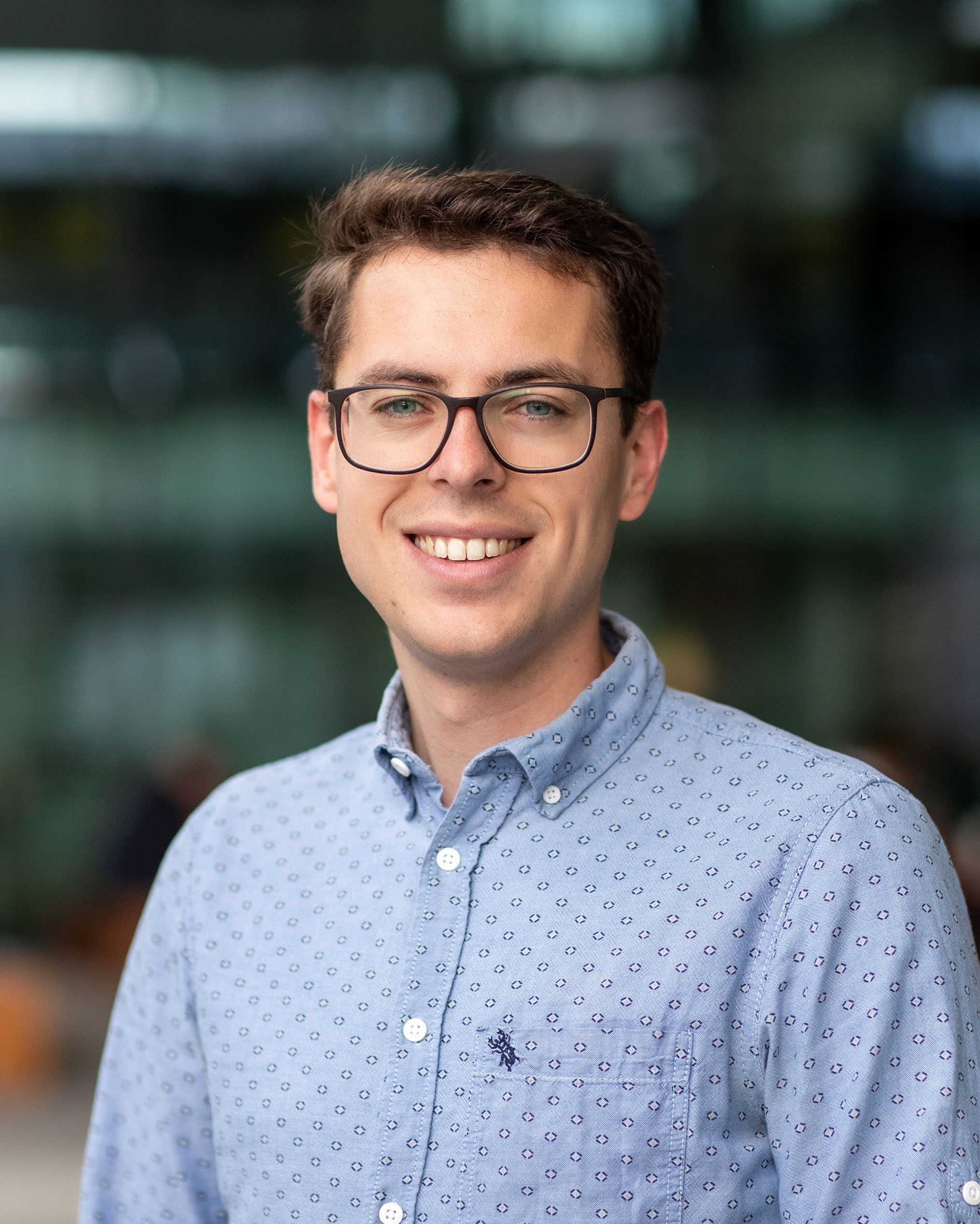}}]{Sven A.N. Nouwens} received the M.Sc. (cum laude) and Ph.D. degrees in mechanical engineering from the Eindhoven University of Technology, Eindhoven, the Netherlands, in 2020 and 2024, respectively. His research interests include numerical methods for model predictive control and the control of partial differential equations.\end{IEEEbiography}\vspace*{-.5cm}
\begin{IEEEbiography}[{\includegraphics[width=1in,height=1.25in,clip,keepaspectratio]{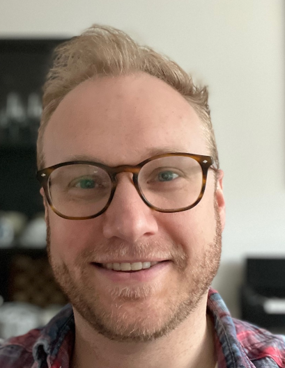}}]{Victor S. Dolk} received the M.Sc. degree (cum laude) in mechanical engineering and the Ph.D. degree (summa cum laude) in control theory from the Eindhoven University of Technology, The Netherlands, in 2013 and 2017, respectively. His research interests include hybrid dynamical systems, networked control systems, event-triggered control, model predictive control and numerical methods for large-scale systems.
\end{IEEEbiography}\vspace*{-.5cm}
\begin{IEEEbiography}[{\includegraphics[width=1in,height=1.25in,clip,keepaspectratio]{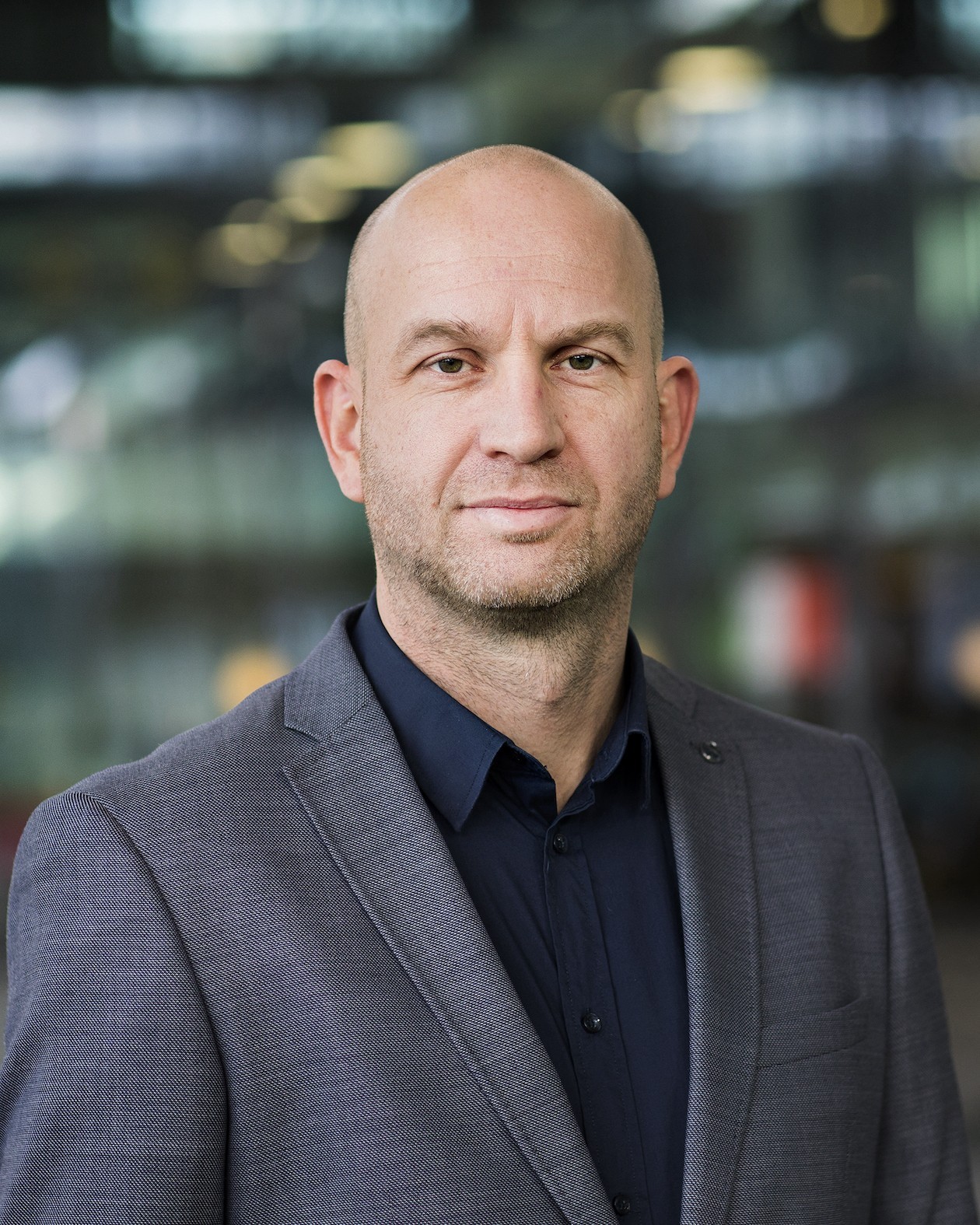}}]{W.P.M.H. (Maurice) Heemels} (Fellow, IEEE) received M.Sc. (mathematics) and Ph.D. (EE, control theory) degrees (summa cum laude) from the Eindhoven University of Technology (TU/e) in 1995 and 1999, respectively. Currently, he is a Full Professor and Vice-Dean at the Department of Mechanical Engineering, TU/e. He held visiting professor positions at ETH, Switzerland (2001), UCSB, USA (2008) and University of Lorraine, France (2020). He is a Fellow of the IEEE and IFAC. He served/s on the editorial boards of Automatica,  Nonlinear Analysis: Hybrid Systems (NAHS), Annual Reviews in Control, and IEEE Transactions on Automatic Control, and is the Editor-in-Chief of NAHS as of 2023. He was the recipient of the 2019 IEEE L-CSS Outstanding Paper Award and the Automatica Paper Prize 2020-2022. His current research includes hybrid and cyber-physical systems, networked and event-triggered control systems and model predictive control.
\end{IEEEbiography}


\end{document}